\documentclass[12pt,reqno]{amsart}

\usepackage{amssymb,amsmath,graphicx,amsfonts,euscript}
\usepackage{color}
\usepackage{appendix}
\usepackage{mathrsfs}

\newtheorem{Theorem}{Theorem}[section]
\newtheorem{Proposition}{Proposition}[section]
\newtheorem{Lemma}{Lemma}[section]
\newtheorem{Corollary}{Corollary}[section]

\newtheorem{Remark}{Remark}[section]
\newcommand{\ls}{\lesssim }

\newcommand{\p}{\partial}
\newcommand{\dtau}{{\rm d} \tau}
\newcommand{\dx}{{\rm d} x}
\newcommand{\dy}{{\rm d} y}

\newcommand{\dki}{{\rm d}\xi}

\numberwithin{equation}{section}

\title[Asymptotic stability of the 2D Boussinesq equations]{Asymptotic stability of the 2D Boussinesq equations without thermal conduction}

\author{Lihua Dong}
\address{Lihua Dong\newline
\noindent Department of Mathematics,\newline
 Nanjing University,\newline
 210093 Nanjing\newline
 People's Republic of China.\newline
e-mail: Donglhmath@smail.nju.edu.cn}
\author{Yongzhong Sun}
\address{Yongzhong Sun\newline
Department of Mathematics, \newline
Nanjing University,\newline
210093 Nanjing\newline
 People's Republic of China.\newline
e-mail: sunyz@nju.edu.cn
}

\begin{document}

\begin{abstract}
 This paper is concerned with the asymptotic stability of certain stationary solution to Boussinesq equations without thermal conduction in the infinite flat strip $\Omega=\mathbb{R}\times (0,1)$. It is shown that the solution starting from initial data close to the stationary solution will converge to it with explicit algebraic rates as time tends to infinity.
 \end{abstract}

\maketitle

\section{Introduction}
As a simplified model, the Boussinesq approximation of Navier-Stokes-Fourier system is widely used in the study of hydrodynamic stability problems, see \cite{PGD1, AJM1, JP1}, among others. Here we consider the 2D Boussinesq system in the absence of thermal conduction. In Eulerian coordinates, it reads as
\begin{equation}
     \label{intp1}
     \left\{
     \begin{aligned}
       &\partial_{t}\mathbf{v} +\mathbf{v}\cdot \nabla \mathbf{v}-\nu\Delta\mathbf{v} +\nabla p = \vartheta\mathbf{e}_{2},\\ 
      &\nabla\cdot\mathbf{v} = 0,\\
 &\partial_{t}\vartheta+\mathbf{v}\cdot \nabla \vartheta = 0.
     \end{aligned}
     \right.
\end{equation}
Here the unknowns  $\mathbf{v} = (v_{1}, v_{2}),\ p$ and $\vartheta$ are the velocity, pressure and temperature of the fluid respectively. The positive constant $\nu$ is the viscosity of the fluid and $\mathbf{e}_{2} = (0,1)$ stands for the direction of buoyancy force.

During the last decades, system (\ref{intp1}) has attracted much interest in the research of mathematical fluid mechanics. The global existence, uniqueness and regularity of smooth solution to (\ref{intp1}) under different settings have been investigated by many authors, see \cite{HA1,DC1,RD2,CRD1,LH1,THSK1,TH1,WH1,NJ1,KW,ML1}, among others. In \cite{CRD1} by C. Doering et al., it is shown that the $L^{2}$ norms of the velocity and its first order derivatives converge to zero as time tends to infinity without any smallness restriction on the initial data. Recently, the authors in \cite{LT1} study the stability of hydrostatic equilibrium to (\ref{intp1}) in the periodic domain and establish decay rates of the velocity under an additional assumption on  the solution. In the recent work \cite{RW2}, the author shows asymptotic behavior and explicit decay rates for solutions to the perturbed system (\ref{Pu4}) in the whole plane $\mathbb{R}^2$. We refer to \cite{LD1} for more detailed introduction and argument on this topic.

In \cite{LD1}, the underlying domain of fluid motion is assumed to be $\Omega=\mathbb{T}\times(0,1)$. Due to spatial confinement in the horizontal direction, in general the temperature can not decay to the stationary one. In the present work we consider the case $\Omega=\mathbb{R}\times (0,1)$. Under such a setting, we show that as time goes to infinity, the solution to (\ref{intp1}) converges to the specific stationary solution with explicit rates provided the initial data is a small perturbation of it. Especially, we obtain the decay rates both for high order derivatives of the velocity and temperature. Moreover, the convergence rates are sharp in the sense that they coincide with that of the linearized equations.  More precisely, we assume that the fluid occupies the two dimensional infinite strip  $\Omega=\mathbb{R}\times(0,1)$ and choose the specific stationary solution $(\vartheta_{s}, \mathbf{v}_s, p_{s})$ to (\ref{intp1}) in $\Omega$ as
\begin{equation}\label{steady0}
\vartheta_{s} = y,\  \mathbf{v}_s = 0,\  p_{s} = \frac{1}{2}y^2,\ y \in [0,1].
\end{equation}
System (\ref{intp1}) is supplemented with the initial conditions
\begin{equation}\label{intpi2}
 \mathbf{v}(\mathbf{x},0) = \mathbf{v}_{0}(\mathbf{x}),\, \vartheta(\mathbf{x},0)=\vartheta_{0}(\mathbf{x})\text{ in }\Omega,
\end{equation}
together with the boundary conditions
\begin{equation}\label{intp2}
 \left(\mathbf{v}\cdot \mathbf{n}\right) (\mathbf{x},t) =0,\  \left(\nabla\times\mathbf{v}\right)(\mathbf{x},t) = 0\text{ on }\partial\Omega,\,t>0,
\end{equation}
where $\mathbf{n}$ is the outward unit normal to $\partial\Omega=\mathbb{R}\times\{y=0,1\}$. By introducing the perturbation
\[
\theta = \vartheta - \vartheta_{s},\  \mathbf{u} = \mathbf{v}-\mathbf{v}_s,\,q= p-p_s,
\]
system (\ref{intp1}) is rewritten as
\begin{equation}\label{Pu4}
     \left\{
     \begin{aligned}
       &\partial_{t}\mathbf{u}+\mathbf{u}\cdot \nabla \mathbf{u} -\nu\Delta\mathbf{u}  + \nabla q = \theta \mathbf{e}_2, \\
       &\nabla\cdot\mathbf{u} = 0,\\
       &\partial_{t}\theta+\mathbf{u}\cdot \nabla \theta = -u_{2},
     \end{aligned}
     \right.
\end{equation}
together with the initial and boundary conditions
\begin{equation}\label{pu42}
 \left\{
     \begin{aligned}
& \mathbf{u}(\mathbf{x},0) = \mathbf{u}_{0}(\mathbf{x}),\, \theta(\mathbf{x},0)=\theta_{0}(\mathbf{x})\text{ in }\Omega,\\
 &\left(\mathbf{u}\cdot \mathbf{n}\right)(\mathbf{x},t) =0,\  \left(\nabla\times\mathbf{u}\right)(\mathbf{x},t) = 0\text{ on }\partial\Omega,\,t>0.
 \end{aligned}
 \right.
\end{equation}
To simplify formulation, we introduce the vorticity $\omega =\partial_{1}u_{2} - \partial_{2}u_{1}$ and the stream function $\psi= (-\Delta)^{-1}\omega$, which solves
\begin{equation}
     \label{P3}
     \left\{
     \begin{aligned}
       & -\Delta\psi = \omega,\\
       &  \psi|_{\partial\Omega} = 0.
     \end{aligned}
     \right.
\end{equation}
Thus  (\ref{Pu4}) is reformulated in the following equations for $(\omega, \theta )$.
\begin{equation}\label{P4}
     \left\{
     \begin{aligned}
       &\partial_{t}\omega-\nu\Delta\omega +\mathbf{u}\cdot \nabla \omega = \partial_{1}\theta, \\
       &\partial_{t}\theta+\mathbf{u}\cdot \nabla \theta = -u_{2},\\
       &\mathbf{u} = (\partial_{2}(-\Delta)^{-1}\omega, -\partial_{1}(-\Delta)^{-1}\omega),\\
     \end{aligned}
     \right.
\end{equation}
The corresponding initial and boundary conditions for $\omega$ are as follows.
\begin{equation}
     \label{P5}
     \left\{
     \begin{aligned}
       & (\omega, \theta)(\mathbf{x},0) = (\omega_{0}, \theta_{0})(\mathbf{x})\text{ in }\Omega,\\
       &  \omega(\mathbf{x},t) = 0\text{ on }\partial\Omega,\,t>0.
     \end{aligned}
     \right.
\end{equation}
We now state the main result of this paper.
\begin{Theorem}\label{Mrt1}
Let $m> 32$ be an integer. Assume that
\begin{eqnarray}
&\omega_{0} \in H^{m}\cap W^{5,1}, \ \partial_{2}^{n}\omega_{0} = 0 \text{ on } \partial\Omega,\text{ for } n = 0, 2, \cdots, 2[(m-1)/2],\label{P61}
\\
&\theta_{0}\in H^{m+1}\cap W^{8,1},\ \partial_{2}^{n}\theta_{0} = 0 \text{ on } \partial\Omega
\text{ for } n = 0, 2, \cdots, 2[m/2]. \label{P6}
\end{eqnarray}
There exists $\epsilon_{0}>0$ depending only on $\nu$ and $m$ such that if
\begin{equation}\label{P9}
\|\theta_{0}\|_{W^{8,1}} + \|\theta_{0}\|_{H^{m+1}} +\|\omega_{0}\|_{W^{5,1}} + \|\omega_{0}\|_{H^{m}}\leq \epsilon_{0},
\end{equation}
then there exists a unique global smooth solution to (\ref{P4})-(\ref{P5}) satisfying
\[
 \|\theta(t)\|_{H^{4}}\lesssim \langle t\rangle^{-\frac{1}{4}},\
\|\mathbf{\omega}(t)\|_{H^{2}}+\|\partial_{1}\theta(t)\|_{H^{2}}\lesssim \langle t\rangle^{-\frac{3}{4}},
\]
\[
\|\partial_{1}\mathbf{\omega}(t)\|_{L^{2}} + \|\partial_{11}\theta(t)\|_{L^{2}}\lesssim \langle t\rangle^{-\frac{5}{4}},
\]
\[
\|\theta(t)\|_{L^{\infty}} + \|\p_{2}\theta(t)\|_{L^{\infty}}\lesssim \langle t\rangle^{-\frac{1}{2}},\ \|\p_{1}\theta(t)\|_{L^{\infty}} +\|\mathbf{\omega}(t)\|_{L^{\infty}}\lesssim \langle t\rangle^{-1}.
\]
\end{Theorem}
\begin{Remark}
\emph{
It is well-known that if a stationary solution $\vartheta_{s}(y)$ satisfies $\vartheta'_{s}(y_0)<0$ for some $y_0\in [0,1]$, which implies that fluid with higher temperature lies below the lower one, then it is unstable-the Rayleigh-Taylor instability happens, see \cite{PB1, CRD1,PGD1}, among others.
}
\end{Remark}
\begin{Remark}
\emph{
As pointed out before, we not only show the asymptotic stability of the stationary solution specified in (\ref{steady0}), but also give the explicit decay rates for high order derivatives of temperature and velocity. Going back to the original initial-boundary value problem (\ref{intp1}), (\ref{intpi2}) and (\ref{intp2}), we have
\[
\|\vartheta(t)-y\|_{H^{4}}\lesssim \langle t\rangle^{-\frac{1}{4}},\,
\|\p_{1}\left(\vartheta(t)-y\right)\|_{H^{2}} +\|v_{1}(t)\|_{H^{3}}\lesssim \langle t\rangle^{-\frac{3}{4}},\
\]
\[
\|\p_{11}\left(\vartheta(t)-y\right)\|_{L^{2}} + \|v_{2}(t)\|_{H^{2}} + \|\partial_{1}v_{1}(t)\|_{H^{1}} \lesssim \langle t\rangle^{-\frac{5}{4}},
\]
\[
\|\vartheta(t)-y\|_{L^{\infty}} + \|\p_{2}\left(\vartheta(t)-y\right)\|_{L^{\infty}}\lesssim \langle t\rangle^{-\frac{1}{2}},\
\|\p_{1}\left(\vartheta(t)-y\right)\|_{L^{\infty}} + \|\mathbf{v}(t)\|_{W^{1,\infty}}\lesssim \langle t\rangle^{-1},
\]
for all $t>0$,
which are consistent with decay rates of the linearized equations (\ref{Eb1more})-(\ref{Eb2}).
}
\end{Remark}

We note that due to the absence of thermal conduction, the mechanism of stabilization in our setting is arising from the action of buoyancy. From the linearized equations (\ref{Eb1more})-(\ref{Eb2}), $\theta$ satisfies
\[
\partial_{tt}\theta-\nu\Delta\partial_{t}\theta - \partial_{1}(-\Delta)^{-1}\partial_{1}\theta = 0,
\]
which exhibits dissipation in the horizontal direction. However, the dissipation is very weak-the time decay rate of $\|\partial_{1}\theta\|_{L^{\infty}}$ is at most $\langle t\rangle^{-1}$, which in turn implies the Lipschitz norm of the velocity has the best decay rate as $\langle t\rangle^{-1}$. Fortunately, this critical decay rate is enough to overcome the difficulty caused by nonlinear terms in the analysis of system (\ref{P4}). Our proof is based on the construction of suitable energy functionals together with a detailed spectral analysis to the linear equation (\ref{Aa1}).

{\bf Notation.}
Throughout this paper, the bold character $\mathbf{x}$ represents the spatial variable $(x,y)\in \mathbb{R}\times (0,1)$. For $m\in \mathbb{N}$, $p\in [1,\infty]$ we denote the inhomogeneous Sobolev space with derivatives up to order $m$ in $L^{p}(\Omega)$ by $W^{m,p}(\Omega)$ equipped with the norm $\|f\|_{W^{m,p}(\Omega)} = \|f\|_{\dot{W}^{m,p}(\Omega)} + \|f\|_{L^p(\Omega)}$, where $\dot{W}^{m,p}(\Omega)$ is the homogeneous Sobolev space with all $m$-th order derivatives belonging to $L^{p}(\Omega)$. Especially, we denote $H^{m}(\Omega)=W^{m,2}(\Omega)$. Notation $\langle\cdot,\cdot\rangle$ is used as the inner product in $L^{2}(\Omega)$. Moreover, the frequency space $\mathbb{R}\times \mathbb{N}$ is denoted as $\widehat{\Omega}$ and $L^{p}(\widehat{\Omega})(1\leq p \leq \infty)$ is the space of all functions $g(\xi,k)$, $(\xi,k)\in \widehat{\Omega}$ with
$\|g\|_{L^{p}(\widehat{\Omega})}^{p} = \sum_{k\in \mathbb{N}}\int_{\mathbb{R}}|g(\xi,k)|^{p}\dki$. For simplicity, we omit the underlying domain by using $L^{p},\,W^{m,p},\,H^{m}$ and $\widehat{L}^{p}$ etc., to denote the corresponding spaces on $\Omega$ and $\widehat\Omega$ respectively. Finally, $A \lesssim B$ means $A\leq CB$ with a generic constant $C$ and $A\sim B$ means $A\lesssim B$ and $B\lesssim A$.
\section{Preliminaries}
In this section, we give some preliminary results that will be used later.

The following estimates on multiplication of two functions in Sobolev spaces and the interpolation inequality are well-known, see \cite{RA1,AJM2}, among others.
\begin{Lemma}\label{PL1}
 Let $m \in \mathbb{N}$.
\begin{itemize}
  \item If $f, g \in H^{m} \cap L^{\infty}$, then
\begin{equation}\label{BL11}
\|fg\|_{H^{m}}\lesssim \|f\|_{H^{m}}\|g\|_{L^{\infty}} + \|f\|_{L^{\infty}}\|g\|_{H^{m}}.
\end{equation}
\item If $f, g \in H^{m}$ with $m >1$, then
\begin{equation}\label{BL12}
\|fg\|_{H^{m}}\lesssim \|f\|_{H^{m}}\|g\|_{H^{m}}.
\end{equation}
\item Let $m_{1}\leq m \leq m_{2}$. If $f\in H^{m_{2}}$, then
\begin{equation}\label{BLT1}
\|f\|_{H^{m}}\lesssim \|f\|_{H^{m_{1}}}^{s}\|f\|_{H^{m_{2}}}^{1-s},\text{ with } m= s m_{1} +(1-s)m_{2},\,0\leq s\leq1.
\end{equation}
\item If $f, g \in H^{m}$, then
\begin{equation}\label{BL13}
\|fg\|_{W^{m,1}}\lesssim \|f\|_{H^{m}}\|g\|_{L^{2}} + \|f\|_{L^{2}}\|g\|_{H^{m}}.
\end{equation}
\item
If $f \in H^{m}\cap W^{1,\infty}$, $g \in H^{m-1} \cap L^{\infty}$, then for $m\ge1$ and any $|\alpha|\leq m$
\begin{equation}\label{BL14}
\|\partial^{\alpha}(fg) - f\partial^{\alpha}g\|_{L^{2}}\lesssim \|\nabla f\|_{L^{\infty}}\|g\|_{H^{m-1}} + \|f\|_{H^{m}}\|g\|_{L^{\infty}}.
\end{equation}
\end{itemize}
\end{Lemma}

For $t>0$, let $\langle t \rangle = \max\{1,t\}$. The following lemma can be verified by a direct calculation.
\begin{Lemma}\label{PL3}
Let $\beta, \gamma>0$ be two constants such that $\beta \leq 1+\gamma$. Then
\begin{equation}\label{BL22}
\int_{0}^{t}\frac{\dtau}{\langle t-\tau\rangle^{\beta}\langle \tau\rangle^{1+\gamma}} \lesssim\langle t\rangle^{-\beta},\,\int_{0}^{t}e^{-( t-\tau)}\langle \tau\rangle^{-\gamma}\dtau \lesssim\langle t\rangle^{-\gamma}.
\end{equation}
\end{Lemma}

We now turn to the Fourier expansion of functions defined in $\Omega=\mathbb{R}\times (0,1)$ with vanishing Dirichlet/Neumann boundary value.
Following \cite{AD1,AD2}, we introduce the functional spaces for $m\in \mathbb{N}$, $p \in [1,\infty]$,
\begin{equation}\label{F1}
 \mathfrak{D}^{m,p} := \{f\in W^{m,p}: \partial_{2}^{n}g|_{\partial\Omega} = 0,\,n = 0, 2, \cdots, 2[(m-1)/2]\},
\end{equation}
\begin{equation}\label{F2}
 \mathfrak{N}^{m,p} := \{f\in W^{m,p}: \partial_{2}^{n}g|_{\partial\Omega} = 0,\,n = 1, 3, \cdots, 2[m/2]-1\},
\end{equation}
and $\mathfrak{D}^{m}=\mathfrak{D}^{m,2}$, $\mathfrak{N}^{m}=\mathfrak{N}^{m,2}$. Here $[m/2]=k$ if $m=2k$ or $2k+1$, $k=1,2,3,\cdots$. The Fourier expansion for a function $f\in \mathfrak{D}^m$ reads as
\begin{equation}\label{pr1}
f(x,y) = \frac{1}{\sqrt{\pi }}\sum_{k=1}^{+\infty}\int_{\mathbb{R}}\widehat{f}_{o}(\xi,k)e^{i x \xi}\dki  \sin k\pi y,
\end{equation}
\[
\widehat{f}_{o}(\xi,k) = \frac{1}{\sqrt{\pi }}\int_{\mathbb{R}}\int_{0}^{1}f(x,y)e^{- i x \xi}\sin k\pi y\dy \dx , \text{ for } (\xi, k)\in \widehat{\Omega},
\]
while for $f\in \mathfrak{N}^m$,
\begin{equation}\label{pr2}
f(x,y) = \frac{1}{\sqrt{\pi }}\sum_{k=0}^{+\infty}\int_{\mathbb{R}}\widehat{f}_{e}(\xi,k)e^{ i x\xi}\dki  \cos k\pi y,
\end{equation}
\[
\widehat{f}_{e}(\xi,k)= \frac{1}{\sqrt{\pi }}\int_{\mathbb{R}}\int_{0}^{1}f(x,y)e^{- i x\xi}\cos k\pi y \dy \dx ,\ \text{ for } (\xi, k)\in \widehat{\Omega}.
\]
We note that $f\in \mathfrak{D}^m $ (or $\mathfrak{N}^{m}$) implies $\partial_2 f \in \mathfrak{N}^{m-1}$ (or $ \mathfrak{D}^{m-1}$) and
\[
\widehat{(\partial_2f)_e}(\xi,k) = -k \pi \widehat{f_o}(\xi,k)\,\left(\widehat{(\partial_2f)_o}(\xi,k) = k\pi \widehat{f_e}(\xi,k)\right).
\]
For notation convenience, we use $\widehat{f}$ to denote $\widehat{f}_{o}$ or $\widehat{f}_{e}$, which makes no confusion once $f \in \mathfrak{D}^m$ or $\mathfrak{N}^m$ is given.
Accordingly,
\begin{equation}\label{eLL2}
\|f\|_{H^m}\sim \|(1+|\cdot|^2)^{\frac{m}{2}}\widehat{f}(\cdot)\|_{\widehat{L}^2},\,f \in \mathfrak{D}^m \text{ or } \mathfrak{N}^m.
\end{equation}
Moreover,
\begin{equation}\label{eLL3}
\|f\|_{L^\infty}\lesssim \|\widehat{f}\|_{\widehat{L}^1},\,\|\widehat{f}\|_{\widehat{L}^\infty}\lesssim \|{f}\|_{{L}^1},\,\|\widehat{f}\|_{\widehat{L}^1} \ls \|f\|_{H^2}\ls \|f\|_{W^{m,1}},\ m \ge 3.
\end{equation}
\begin{equation}\label{eLL4}
\|(1+|\cdot|^2)^{\frac{m}{2}}\widehat{f}\|_{\widehat{L}^\infty}\lesssim \|{f}\|_{{W}^{m,1}},\, f \in \mathfrak{D}^{m,1} \text{ or } \mathfrak{N}^{m,1}.
\end{equation}
Note that we use $\|\widehat{f}(\cdot)\|_{\widehat{H}^m}$ to denote $\|(1+|\cdot|^2)^{\frac{m}{2}}\widehat{f}(\cdot)\|_{\widehat{L}^2}$. In the following context, these facts will be used frequently without being referred.

Finally, we consider the following Dirichlet problem for Laplace equation in $\Omega$.
\begin{equation}
     \label{Poise}
     \left\{
     \begin{aligned}
       & -\Delta\varphi= f,  \\
       &  \varphi|_{\partial\Omega} = 0.
     \end{aligned}
     \right.
\end{equation}
For $f\in L^2$, this boundary value problem is uniquely solvable in $H_0^1\cap H^2$. We denote this unique solution $ \varphi = (-\Delta)^{-1} f$.
\begin{Lemma}\label{Pois1}
Assume that $f \in H^{m}$ for $m \ge 0$.
Then
\begin{equation}\label{eLL1}
\|\varphi\|_{H^{m+2}}\lesssim \|f\|_{H^{m}}.
\end{equation}
\end{Lemma}
\begin{Remark}
\emph{
Inequality (\ref{eLL1}) is nothing but the standard elliptic estimate. For $f\in \mathfrak{D}^m$, (\ref{eLL1}) immediately follows from (\ref{eLL2}) and
\[
\widehat{\varphi}_o(\xi,k) = \frac{1}{\xi^{2}+\pi^{2}k^{2}}\widehat{f}_o(\xi,k),\, \xi\in \mathbb{R},\, k \ge 1.
\]
}
\end{Remark}

By using the stream function $\psi$ in (\ref{P3}), together with the help of Lemma \ref{Pois1}, we have
\begin{Corollary}\label{corollary1}
Assume that $\mathbf{u}\in H^1$ satisfying $\nabla\cdot \mathbf{u} =0$ in $\Omega$ and $\mathbf{u}\cdot \mathbf{n} = 0$ on $\partial\Omega$. If
$\omega=\partial_1 u_2 -\partial_2 u_1 \in H^{m}$, $m \in \mathbb{N}$, then $\mathbf{u}\in H^{m+1}$ and
\begin{equation}\label{FC1}
\|\mathbf{u}\|_{H^{m+1}}\lesssim \|\omega\|_{H^{m}}.
\end{equation}
\end{Corollary}
\section{Decay estimates for solutions to the linearized equations}
The linearized equations of (\ref{P4}) is
\begin{equation}\label{Eb1more}
     \left\{
     \begin{aligned}
       &\partial_{t}w-\nu\Delta w = \partial_{1}\phi, \\
       &\partial_{t}\phi + u_{2} = 0, \\
       & u_{2} = -\partial_{1}(-\Delta)^{-1} w,
     \end{aligned}
     \right.
\end{equation}
together with the initial and boundary conditions
\begin{equation}
     \label{Eb2}
     \left\{
     \begin{aligned}
       & (w, \phi)(\mathbf{x},0) = (w_{0}, \phi_{0})(\mathbf{x})\text{ in }\Omega, \\
       &  w(\mathbf{x},t) = 0 \text{ on }\partial\Omega,\,t>0.
     \end{aligned}
     \right.
\end{equation}
In this section we give explicit decay estimates of the solutions to (\ref{Eb1more})-(\ref{Eb2}).
\subsection{Decay rates for solutions to the linearized equation of temperature}
We decouple system $(\ref{Eb1more})$ by taking the time derivative on $(\ref{Eb1more})_{2}$ and using the remaining equations of $(\ref{Eb1more})$ to find
\begin{equation}\label{Eb1}
    \partial_{tt}\phi-\nu\Delta\partial_{t}\phi - \partial_{1}(-\Delta)^{-1}\partial_{1}\phi = 0.
\end{equation}
We consider the following inhomogeneous equation in $\Omega$,
\begin{equation}\label{Aa1}
   \partial_{tt}\phi-\nu\Delta\partial_{t}\phi - \partial_{1}(-\Delta)^{-1}\partial_{1}\phi =F
\end{equation}
with the initial and boundary conditions
\begin{equation}\label{Aam1}
\left\{
     \begin{aligned}
       & \phi(\mathbf{x},0)=\phi_{0}(\mathbf{x}), \partial_{t}\phi(\mathbf{x},0)=\phi_{1}(\mathbf{x})\text{ in }\Omega,\\
       &  \phi(\mathbf{x},t)=0 \text{ on }\partial\Omega,\,t>0.
     \end{aligned}
     \right.
\end{equation}
To solve the above initial-boundary problem, for $t>0$ we define operators $\mathcal{L}_{1}(t)$ and $\mathcal{L}_{2}(t)$ as follows.
\begin{equation}\label{Aa2}
  \widehat{\mathcal{L}_{1}(t)f}(\xi,k) = \frac{1}{2}\left(e^{-\frac{1}{2}\left(\nu(\xi^{2}+ \pi^{2}k^{2})+\sigma\right)t}+e^{-\frac{1}{2}\left(\nu(\xi^{2}+ \pi^{2} k^{2})-\sigma\right)t}\right)\widehat{f}(\xi,k),
\end{equation}
and
\begin{equation}\label{Aa3}
\widehat{\mathcal{L}_{2}(t)f}(\xi,k) = \frac{1}{\sigma}\left(e^{-\frac{1}{2}\left(\nu(\xi^{2}+ \pi^{2}k^{2})-\sigma\right)t}-e^{-\frac{1}{2}\left(\nu(\xi^{2}+ \pi^{2}k^{2})+\sigma\right)t}\right)\widehat{f}(\xi,k).
\end{equation}
Here $\sigma= \sqrt{\nu^{2}(\xi^{2}+ \pi^{2}k^{2})^{2}-\frac{4\xi^{2}}{\xi^{2}+ \pi^{2}k^{2}}}$ and $(\xi,k)\in\widehat{\Omega}$.
\begin{Lemma}\label{AL1}
Assume that $\phi_0\in H^2\cap H_0^1$, $\phi_1\in L^2$ and $F\in L_{loc}^1(0, \infty; L^2)$. Then the solution to (\ref{Aa1}) is given by
\[
\phi(x,y,t)= \mathcal{L}_{1}(t)\phi_{0}(x,y) + \mathcal{L}_{2}(t)\left(\frac{\nu}{2}(-\Delta)\phi_{0}(x,y) +\phi_{1}(x,y)\right)
\]
\begin{equation}\label{sol}
+ \int_{0}^{t}\mathcal{L}_{2}(t-\tau)F(x,y,\tau)\dtau.
\end{equation}
\end{Lemma}
{\bf Proof.}
Performing the Fourier transform to (\ref{Aa1}),
\begin{equation}\label{Aa4}
   \partial_{tt}\widehat{\phi}(\xi,k,t)+ \nu\left(\xi^{2}+ \pi^{2}k^{2}\right)\partial_{t}\widehat{\phi}(\xi,k,t)+\frac{\xi^{2}}{\xi^{2}+ \pi^{2}k^{2}}\widehat{\phi}(\xi,k,t)=\widehat{F}(\xi,k,t), \, (\xi,k)\in\widehat{\Omega}.
\end{equation}
Note that
\[
\left(\partial_{tt} + \nu\left(\xi^{2}+ \pi^{2}k^{2}\right)\partial_{t} +\frac{\xi^{2}}{\xi^{2}+ \pi^{2}k^{2}}\right)\widehat{\phi}(\xi,k,t)
\]
\[
 = \left[ \left(\partial_{t}+\frac{\nu(\xi^{2}+ \pi^{2}k^{2})}{2}\right)^{2}-\frac{\nu^{2}(\xi^{2}+ \pi^{2}k^{2})^{2}}{4}+\frac{\xi^{2}}{\xi^{2}+ \pi^{2}k^{2}}\right]\widehat{\phi}(\xi,k,t)
\]
\[
=\left(\partial_{t}+\frac{\nu(\xi^{2}+\pi^{2}k^{2})}{2}+\frac{1}{2}\sqrt{\nu^{2}(\xi^{2}+ \pi^{2}k^{2})^{2}-\frac{4\xi^{2}}{\xi^{2}+ \pi^{2}k^{2}}}\right)
\]
\[
\cdot\left(\partial_{t}+\frac{\nu(\xi^{2}+ \pi^{2}k^{2})}{2}-\frac{1}{2}\sqrt{\nu^{2}(\xi^{2}+ \pi^{2}k^{2})^{2}-\frac{4\xi^{2}}{\xi^{2}+ \pi^{2}k^{2}}}\right)\widehat{\phi}(\xi,k,t).
\]
Let
\begin{equation}\label{EA+}
  \psi_{+}(\xi,k,t)=\left(\partial_{t}+\frac{\nu(\xi^{2}+ \pi^{2}k^{2})}{2}+\frac{\sigma}{2}\right)\widehat{\phi}(\xi,k,t)
\end{equation}
and
\begin{equation}\label{EA++}
\psi_{-}(\xi,k,t)=\left(\partial_{t}+\frac{\nu(\xi^{2}+ \pi^{2}k^{2})}{2}-\frac{\sigma}{2}\right)\widehat{\phi}(\xi,k,t).
\end{equation}
Then
\begin{equation}\label{2+1}
\left(\partial_{t}+\frac{\nu(\xi^{2}+ \pi^{2}k^{2})}{2}-\frac{\sigma}{2}\right)\psi_{+}=\widehat{F},\
\left(\partial_{t}+\frac{\nu(\xi^{2}+ \pi^{2}k^{2})}{2}+\frac{\sigma}{2}\right)\psi_{-}=\widehat{F}.
\end{equation}
By Duhamel's principle,
\begin{equation}\label{1+3}
\psi_{+}(\xi,k,t) = e^{\left(-\frac{\nu(\xi^{2}+ \pi^{2}k^{2})}{2}+\frac{\sigma}{2}\right)t}\psi_{+}(\xi,k,0)
+\int_{0}^{t}e^{\left(-\frac{\nu(\xi^{2}+ \pi^{2}k^{2})}{2}+\frac{\sigma}{2}\right)(t-\tau)}\widehat{F}(\xi,k,\tau)\dtau,
\end{equation}
\begin{equation}\label{1+2}
\psi_{-}(\xi,k,t) = e^{\left(-\frac{\nu(\xi^{2}+ \pi^{2}k^{2})}{2}-\frac{\sigma}{2}\right)t}\psi_{-}( \xi,k,0)
+\int_{0}^{t}e^{\left(-\frac{\nu(\xi^{2}+ \pi^{2}k^{2})}{2}-\frac{\sigma}{2}\right)(t-\tau)}\widehat{F}(\xi,k,\tau)\dtau.
\end{equation}
Moreover, from (\ref{EA+}) and (\ref{EA++}), it follows that
\begin{equation}\label{1+5}
\psi_{+}(\xi,k,0)=\widehat{\phi}_{1}(\xi, k)+\frac{\nu(\xi^{2}+ \pi^{2}k^{2})}{2}\widehat{\phi}_{0}(\xi, k)+\frac{\sigma}{2}\widehat{\phi}_{0}(\xi, k),
\end{equation}
\begin{equation}\label{1+4}
 \psi_{-}(\xi,k,0)=\widehat{\phi}_{1}(\xi, k)+\frac{\nu(\xi^{2}+ \pi^{2}k^{2})}{2}\widehat{\phi}_{0}(\xi, k)-\frac{\sigma}{2}\widehat{\phi}_{0}(\xi, k).
\end{equation}
Substituting (\ref{1+5}) and (\ref{1+4}) into (\ref{1+3}) and (\ref{1+2}) respectively, together with the fact that $\widehat{\phi} = \frac{\psi_{+}-\psi_{-}}{\sigma}$,
\[
\widehat{\phi} =   \frac{1}{2}\left(e^{-\frac{1}{2}\left(\nu(\xi^{2}+ \pi^{2}k^{2})+\sigma\right)t}+e^{-\frac{1}{2}\left(\nu(\xi^{2}+ \pi^{2}k^{2})-\sigma\right)t}\right)\widehat{\phi}_{0}
\]
\[
+ \frac{1}{\sigma}\left(e^{-\frac{1}{2}\left(\nu(\xi^{2}+ \pi^{2}k^{2})-\sigma\right)t}-e^{-\frac{1}{2}\left(\nu(\xi^{2}+ \pi^{2}k^{2})+\sigma\right)t}\right)\left(\widehat{\varphi}_{0}+\frac{\nu}{2}(\xi^{2}+ \pi^{2}k^{2})\widehat{\phi}_{0}\right)
\]
\[ +\int_{0}^{t}\frac{1}{\sigma}\left(e^{-\frac{1}{2}\left(\nu(\xi^{2}+ \pi^{2}k^{2})-\sigma\right)(t-\tau)}-e^{-\frac{1}{2}\left(\nu(\xi^{2}+ \pi^{2}k^{2})+\sigma\right)(t-\tau)}\right)\widehat{F}(\tau)\dtau,
\]
which is nothing but (\ref{sol}).
\hfill$\square$

The next lemma is inspired by \cite{TE1}.
\begin{Lemma}\label{AL0}
Let $g\in \mathfrak{D}^{6,1}$ and
\[
   \widehat{G}(\xi,k,t) = e^{-\frac{\xi^{2}}{\left(\xi^{2}+\pi^{2}k^{2}\right)^{2}}t}\widehat{g}(\xi,k), \,(\xi,k) \in \widehat{\Omega}.
\]
Then
\begin{equation}\label{P10}
\|\widehat{G}(t)\|_{\widehat{L}^{1}}\lesssim \langle t \rangle^{-\frac{1}{2}}\|g\|_{W^{4,1}},\  \| \widehat{\p_1G}(t)\|_{\widehat{L}^{1}}\lesssim \langle t \rangle^{-1}\|g\|_{W^{6,1}},
\end{equation}
\begin{equation}\label{P101}
\|\widehat{G}(t)\|_{\widehat{L}^{2}}\lesssim \langle t \rangle^{-\frac{1}{4}}\|g\|_{W^{2,1}},\  \| \widehat{\p_1G}(t)\|_{\widehat{L}^{2}}\lesssim \langle t \rangle^{-\frac{3}{4}}\|g\|_{W^{4,1}}.
\end{equation}
\end{Lemma}
{\bf Proof.}
A straightforward calculation shows that
\[
\|\widehat{G}(t)\|_{\widehat{L}^{1}}
\lesssim \sum_{k=1}^{+\infty}\int_{\mathbb{R}}e^{-\frac{\xi^{2}}{\left(\xi^{2}+ \pi^{2}k^{2}\right)^{2}}t}|\widehat{g}(\xi,k)|\dki
\]
\[
\lesssim \|\left(\xi^{2}+\pi^{2}k^{2}\right)^{2}\widehat{g}(\xi,k)\|_{\widehat{L}^{\infty}}\sum_{k=1}^{+\infty}\int_{\mathbb{R}}e^{-\frac{\xi^{2}}{\left(\xi^{2}+ \pi^{2}k^{2}\right)^{2}}t}\left(\xi^{2}+\pi^{2} k^{2}\right)^{-2}\dki
\]
\[
\lesssim \int_{\pi}^{+\infty} \int_{\mathbb{R}}e^{-\frac{\xi^{2}}{\left(\xi^{2}+\eta^{2}\right)^{2}}t}\left(\xi^{2}+\eta^{2}\right)^{-2}\dki {\rm d}\eta\|g\|_{W^{4,1}}.
\]
Using polar coordinates,
\[
\int_{\pi}^{+\infty}\int_{\mathbb{R}}e^{-\frac{\xi^{2}}{\left(\xi^{2}+\eta^{2}\right)^{2}}t}\left(\xi^{2}+\eta^{2}\right)^{-2}\dki {\rm d}\eta
\lesssim \int_{\pi}^{+\infty}\int_{0}^{2\pi}e^{-\frac{\cos^{2}\beta}{r^{2}}t}r^{-3}{\rm d}\beta {\rm d}r.
\]
We decompose the last integral into
\[
\sum_{j=1}^{8}\int_{\pi}^{+\infty}\int_{\frac{(j-1)\pi}{4}}^{\frac{j \pi}{4}}e^{-\frac{\cos^{2}\beta}{r^{2}}t}r^{-3}{\rm d}\beta {\rm d}r = \sum_{j=1}^{8} M_j.
\]
It is enough to consider, say
\[
M_{2}=\int_{\pi}^{+\infty}\int_{\frac{\pi}{4}}^{\frac{\pi}{2}}e^{-\frac{\cos^{2}\beta}{r^{2}}t}r^{-3}{\rm d}\beta {\rm d}r,\,
M_{4}=\int_{\pi}^{+\infty}\int_{\frac{3\pi}{4}}^{\pi}e^{-\frac{\cos^{2}\beta}{r^{2}}t}r^{-3}{\rm d}\beta {\rm d}r.
\]
It is obvious that $M_2 \ls 1$ and $M_4 \ls 1$. By using change of variable $z= \frac{\sqrt{t}}{r}\cos \beta $,
\[
M_{2}\lesssim t^{-\frac{1}{2}}\int_{\pi}^{+\infty}\int_{0}^{\frac{\sqrt{t}}{r}}e^{-z^{2}}r^{-2}\frac{1}{\sqrt{1-\frac{r^{2}z^{2}}{t}}}{\rm d}z{\rm d}r
\]
\[
\lesssim t^{-\frac{1}{2}}\int_{\pi}^{+\infty}\frac{1}{r^{2}}dr\int_{-\infty}^{+\infty}e^{-z^{2}}dz
\lesssim t^{-\frac{1}{2}}.
\]
Hence,
\[
M_2 \ls \langle t\rangle^{-\frac{1}{2}}.
\]
For $M_4$, first note that
\[
\frac{\sqrt{2}}{2}\leq |\cos \beta| \leq 1 \text{ for all }\beta \in \left[\frac{3\pi}{4},\pi\right].
\]
Then
\[
M_{4}
\lesssim \int_{\pi}^{+\infty}\int_{\frac{3\pi}{4}}^{\pi}e^{-\frac{t}{2r^{2}}}r^{-3}{\rm d}\beta {\rm d}r
= \int_{\pi}^{+\infty}\int_{\frac{3\pi}{4}}^{\pi}e^{-\frac{t}{2r^{2}}}\frac{t^{\frac{1}{2}}}{r}\frac{r}{t^{\frac{1}{2}}}r^{-3}{\rm d}\beta {\rm d}r
\]
\[
\lesssim t^{-\frac{1}{2}}\int_{\pi}^{+\infty}\frac{1}{r^{2}}{\rm d}r
\lesssim t^{-\frac{1}{2}}.
\]
Thus, the first estimate in (\ref{P10}) is obtained.
In a similar way we have
\[
\|\widehat{\p_1G}(t)\|_{\widehat{L}^{1}}
\lesssim \sum_{k=1}^{+\infty}\int_{\mathbb{R}}e^{-\frac{\xi^{2}}{\left(\xi^{2}+ \pi^{2}k^{2}\right)^{2}}t}|\xi\widehat{g}(\xi,k)|\dki
\]
\[
\lesssim t^{-\frac{1}{2}}\sum_{k=1}^{+\infty}\int_{\mathbb{R}}e^{-\frac{\xi^{2}}{\left(\xi^{2}+ \pi^{2}k^{2}\right)^{2}}t}\frac{|\xi|}{\xi^{2}+ \pi^{2}k^{2}}t^{\frac{1}{2}}\left(\xi^{2}+ \pi^{2}k^{2}\right)|\widehat{g}(\xi,k)|\dki
\]
\[
\lesssim t^{-\frac{1}{2}}\sum_{k=1}^{+\infty}\int_{\mathbb{R}}e^{-\frac{\xi^{2}}{2\left(\xi^{2}+ \pi^{2}k^{2}\right)^{2}}t}\left(\xi^{2}+ \pi^{2}k^{2}\right)|\widehat{g}(\xi,k)|\dki
 \lesssim \langle t \rangle^{-1}\|g\|_{W^{6,1}}.
\]
Thus  (\ref{P10}) has been verified. Moreover,
\[
\|\widehat{G}(t)\|_{\widehat{L}^{2}}^{2}
\lesssim \sum_{k=1}^{+\infty}\int_{\mathbb{R}}e^{-\frac{2\xi^{2}}{\left(\xi^{2}+ \pi^{2}k^{2}\right)^{2}}t}|\widehat{g}(\xi,k)|^{2}\dki
\]
\[
\lesssim \|\left(\xi^{2}+\pi^{2}k^{2}\right)^{2}\widehat{g}^{2}(\xi,k)\|_{\widehat{L}^{\infty}}\sum_{k=1}^{+\infty}\int_{\mathbb{R}}e^{-\frac{2\xi^{2}}{\left(\xi^{2}+ \pi^{2}k^{2}\right)^{2}}t}\left(\xi^{2}+\pi^{2} k^{2}\right)^{-2}\dki
\]
\[
\lesssim \int_{\pi}^{+\infty} \int_{\mathbb{R}}e^{-\frac{2\xi^{2}}{\left(\xi^{2}+\eta^{2}\right)^{2}}t}\left(\xi^{2}+\eta^{2}\right)^{-2}\dki {\rm d}\eta\|g\|_{W^{2,1}}^{2}
\ls \langle t \rangle^{-\frac{1}{2}}\|g\|_{W^{2,1}}^{2}.
\]
Hence the first estimate in (\ref{P101}) has been proved. The proof of $(\ref{P101})_2$ is similar.
\hfill$\square$
\begin{Remark}\label{Rema1}
\emph{
In the proof of Lemma {\ref{AL0}}, we in fact use $\|g\|_{W^{m,1}}$ to control
\[
\|(1+|\cdot|^2)^{\frac{m}{2}}\widehat{g}\|_{\widehat{L}^{\infty}}
\]
according to (\ref{eLL4}).
}
\end{Remark}

We proceed to give decay estimates of $\mathcal{L}_{1}(t)$, $\mathcal{L}_{2}(t)$.
\begin{Lemma}\label{AL2}
Let $f\in \mathfrak{D}^{6,1}$. Then for any fixed $\nu >0$,
\begin{equation}\label{All1}
\left\|\widehat{\mathcal{L}_{1}(t)f}\right\|_{\widehat{L}^{1}} + \left\|\widehat{\mathcal{L}_{2}(t)f}\right\|_{\widehat{L}^{1}}\lesssim \langle t \rangle^{-\frac{1}{2}}\|f\|_{W^{4,1}},
\end{equation}
\begin{equation}\label{All1+}
\left\|\widehat{\partial_2\mathcal{L}_{1}(t)f}\right\|_{\widehat{L}^{1}} + \left\|\widehat{\partial_2\mathcal{L}_{2}(t)f}\right\|_{\widehat{L}^{1}}\lesssim \langle t \rangle^{-\frac{1}{2}}\|f\|_{W^{5,1}},
\end{equation}
\begin{equation}\label{All2}
\left\|\widehat{\partial_{t}\mathcal{L}_{1}(t)f}\right\|_{\widehat{L}^{1}} + \left\|\widehat{\partial_{t}\mathcal{L}_{2}(t)f}\right\|_{\widehat{L}^{1}}\lesssim \langle t \rangle^{-\frac{3}{2}}\|f\|_{W^{5,1}},
\end{equation}
\begin{equation}\label{All3}
\left\|\widehat{\p_1\mathcal{L}_{1}(t)f}\right\|_{\widehat{L}^{1}} + \left\|\widehat{\p_1\mathcal{L}_{2}(t)f}\right\|_{\widehat{L}^{1}}\lesssim \langle t \rangle^{-1}\|f\|_{W^{6,1}},
\end{equation}
\begin{equation}\label{All4}
\left\|\widehat{\mathcal{L}_{1}(t)f}\right\|_{\widehat{L}^{2}} + \left\|\widehat{\mathcal{L}_{2}(t)f}\right\|_{\widehat{L}^{2}}\lesssim \langle t \rangle^{-\frac{1}{4}}\|f\|_{W^{2,1}},
\end{equation}
\begin{equation}\label{All4+}
\left\|\widehat{\partial_{2}\mathcal{L}_{1}(t)f}\right\|_{\widehat{L}^{2}} + \left\|\widehat{\partial_{2}\mathcal{L}_{2}(t)f}\right\|_{\widehat{L}^{2}}\lesssim \langle t \rangle^{-\frac{1}{4}}\|f\|_{W^{3,1}},
\end{equation}
\begin{equation}\label{All5}
\left\|\widehat{\partial_{t}\mathcal{L}_{1}(t)f}\right\|_{\widehat{L}^{2}} + \left\|\widehat{\partial_{t}\mathcal{L}_{2}(t)f}\right\|_{\widehat{L}^{2}}\lesssim \langle t \rangle^{-\frac{5}{4}}\|f\|_{W^{3,1}},
\end{equation}
\begin{equation}\label{All7}
\left\|\widehat{\p_1\mathcal{L}_{1}(t)f}\right\|_{\widehat{L}^{2}} + \left\|\widehat{\p_1\mathcal{L}_{2}(t)f}\right\|_{\widehat{L}^{2}}\lesssim \langle t \rangle^{-\frac{3}{4}}\|f\|_{W^{4,1}},
\end{equation}
\begin{equation}\label{All6}
\left\|\widehat{\p_{11}\mathcal{L}_{1}(t)f}\right\|_{\widehat{L}^{2}} + \left\|\widehat{\p_{11}\mathcal{L}_{2}(t)f}\right\|_{\widehat{L}^{2}}\lesssim \langle t \rangle^{-\frac{5}{4}}\|f\|_{W^{6,1}}.
\end{equation}
\end{Lemma}
{\bf Proof.}
Let
\[
\lambda_{+}= \frac{-\nu(\xi^{2}+\pi^{2}k^{2})+\sigma}{2},~\lambda_{-}= \frac{-\nu(\xi^{2}+\pi^{2}k^{2})-\sigma}{2}.
\]
A straightforward calculation shows that
\[
\widehat{\mathcal{L}_{1}(t)} = \frac{1}{2}\left(e^{\lambda_{+}t}+e^{\lambda_{-}t}\right),\
\widehat{\mathcal{L}_{2}(t)}= \frac{1}{\lambda_{+}-\lambda_{-}}\left(e^{\lambda_{+}t}-e^{\lambda_{-}t}\right),
\]
\[
\widehat{\partial_{t}\mathcal{L}_{1}(t)} = \frac{1}{2}\left(\lambda_{+}e^{\lambda_{+}t}+\lambda_{-}e^{\lambda_{-}t}\right),\
\widehat{\partial_{t}\mathcal{L}_{2}(t)}= \frac{1}{\lambda_{+}-\lambda_{-}}\left(\lambda_{+}e^{\lambda_{+}t}-\lambda_{-}e^{\lambda_{-}t}\right),
\]
where $\lambda_{+}-\lambda_{-} = \sigma$.

Note  that
\[
\nu_{*}^{2} :=\sup_{\xi \in \mathbb{R},k\geq 1}\frac{4\xi^{2}}{(\xi^{2}+\pi^{2}k^{2})^{3}} \in (0,1).
\]
We discuss decay estimates of $\mathcal{L}_{1}(t)$, $\mathcal{L}_{2}(t)$ in two cases: $\nu \geq\nu_{*}$ and $0<\nu <\nu_{*}$.

First, for $\nu \geq \nu_{*}$,
\[
\lambda_{+} = -\frac{1}{2}\left(\nu(\xi^{2}+\pi^{2}k^{2})-\sigma\right) = -\frac{1}{2}\left(\nu(\xi^{2}+\pi^{2}k^{2})-\sqrt{\nu^{2}(\xi^{2}+\pi^{2}k^{2})^{2}-\frac{4\xi^{2}}{\xi^{2}+\pi^{2}k^{2}}}\right)
\]
\[
= -\frac{2\xi^{2}}{\nu(\xi^{2}+\pi^{2}k^{2})^{2}}\frac{1}{1+\sqrt{1-\frac{4\xi^{2}}{\nu^{2}(\xi^{2}+\pi^{2}k^{2})^{3}}}}.
\]
Hence
\[
-\frac{2\xi^{2}}{\nu(\xi^{2}+\pi^{2}k^{2})^{2}}\leq \lambda_{+} \leq -\frac{\xi^{2}}{\nu(\xi^{2}+\pi^{2}k^{2})^{2}},\
-\nu(\xi^{2}+\pi^{2}k^{2})\leq\lambda_{-} \leq -\frac{\nu}{2}(\xi^{2}+\pi^{2}k^{2}).
\]
Moreover, there exists $\lambda_{0}>0$ such that
\[
\lambda_{0}\leq\lambda_{+} -\lambda_{-} < \nu(\xi^{2}+\pi^{2}k^{2}).
\]
For $k\geq 1$,
\[
0<\widehat{\mathcal{L}_{1}(t)}\lesssim e^{-\frac{\xi^{2}}{\nu(\xi^{2}+\pi^{2}k^{2})^{2}}t},
\]
\[
0\leq \widehat{\mathcal{L}_{2}(t)}\lesssim \frac{1}{\lambda_{0}}e^{-\frac{\xi^{2}}{\nu(\xi^{2}+\pi^{2}k^{2})^{2}}t} + \frac{1}{\lambda_{0}}e^{-\frac{\nu(\xi^{2}+\pi^{2}k^{2})}{2}t}
\lesssim e^{-\frac{\xi^{2}}{\nu(\xi^{2}+\pi^{2}k^{2})^{2}}t},
\]
\[
\left|\widehat{\partial_{t}\mathcal{L}_{1}(t)}\right|\ls
\frac{\xi^{2}}{\nu(\xi^{2}+\pi^{2}k^{2})^{2}}e^{-\frac{\xi^{2}}{\nu(\xi^{2}+\pi^{2}k^{2})^{2}}t}
+\nu (\xi^{2}+\pi^{2}k^{2})e^{-\frac{\nu(\xi^{2}+\pi^{2}k^{2})}{2}t},
\]
\[
\left|\widehat{\partial_{t}\mathcal{L}_{2}(t)}\right|\ls
\frac{\xi^{2}}{\nu(\xi^{2}+\pi^{2}k^{2})^{2}}e^{-\frac{\xi^{2}}{\nu(\xi^{2}+\pi^{2}k^{2})^{2}}t}
+\nu(\xi^{2}+\pi^{2}k^{2})e^{-\frac{\nu(\xi^{2}+\pi^{2}k^{2})}{2}t}.
\]
We use (\ref{eLL3}) and Lemma \ref{AL0} to obtain
\[
\left\|\widehat{\mathcal{L}_{1}(t)f}\right\|_{\widehat{L}^{1}}
+\left\|\widehat{\mathcal{L}_{2}(t)f}\right\|_{\widehat{L}^{1}}
\]
\begin{equation}\label{tg1}
\ls \sum_{k=1}^{\infty}\int_{\mathbb{R}}e^{-\frac{\xi^{2}}{\nu(\xi^{2}+\pi^{2}k^{2})^{2}}t}|\widehat{f}(\xi,k)|\dki
\ls \langle t \rangle^{-\frac{1}{2}}\|f\|_{W^{4,1}},
\end{equation}
\[
\left\|\widehat{\partial_2\mathcal{L}_{1}(t)f}\right\|_{\widehat{L}^{1}}
+\left\|\widehat{\partial_2\mathcal{L}_{2}(t)f}\right\|_{\widehat{L}^{1}}
\]
\begin{equation}\label{tg2}
\ls \sum_{k=1}^{\infty}\int_{\mathbb{R}}e^{-\frac{\xi^{2}}{\nu(\xi^{2}+\pi^{2}k^{2})^{2}}t}|k\widehat{f}(\xi,k)|\dki
\ls \langle t \rangle^{-\frac{1}{2}}\|f\|_{W^{5,1}},
\end{equation}
\[
\left\|\widehat{\partial_{t}\mathcal{L}_{1}(t)f}\right\|_{\widehat{L}^{1}}+\left\|\widehat{\partial_{t}\mathcal{L}_{2}(t)f}\right\|_{\widehat{L}^{1}}
\]
\[
\ls \sum_{k=1}^{\infty}\int_{\mathbb{R}}\frac{\xi^{2}}{\nu(\xi^{2}+\pi^{2}k^{2})^{2}}e^{-\frac{\xi^{2}}{\nu(\xi^{2}+\pi^2k^{2})^{2}}t}|\widehat{f}(\xi,k)|\dki
\]
\[
+ \sum_{k=1}^{\infty}\int_{\mathbb{R}}\nu(\xi^{2}+\pi^{2}k^{2})e^{-\frac{\nu(\xi^{2}+\pi^{2}k^{2})}{2}t}|\widehat{f}(\xi,k)|\dki
\]
\begin{equation}\label{tg4}
\ls \langle t \rangle^{-\frac{3}{2}}\|f\|_{W^{4,1}} + e^{-\frac{\nu}{2}t}\|(\xi^{2}+\pi^{2}k^{2})\widehat{f}\|_{\widehat{L}^{1}}
\ls \langle t \rangle^{-\frac{3}{2}}\|f\|_{W^{5,1}},
\end{equation}
\[
\left\|\widehat{\p_1\mathcal{L}_{1}(t)f}\right\|_{\widehat{L}^{1}}
+\left\|\widehat{\p_1\mathcal{L}_{2}(t)f}\right\|_{\widehat{L}^{1}}
\]
\begin{equation}\label{tg5}
\ls \sum_{k=1}^{\infty}\int_{\mathbb{R}}|\xi|e^{-\frac{\xi^{2}}{\nu(\xi^{2}+\pi^{2}k^{2})^{2}}t}|\widehat{f}(\xi,k)|\dki
\ls \langle t \rangle^{-1}\|f\|_{W^{6,1}}.
\end{equation}
Other estimates in (\ref{All4})-(\ref{All6}) can be obtained in a similar way.

Next, for $0<\nu <\nu_{*}$,
$
\sigma= \sqrt{\nu^2(\xi^{2}+ \pi^{2}k^{2})^{2}-\frac{4\xi^{2}}{\xi^{2}+ \pi^{2}k^{2}}}
$
is not necessary a real number for any $(\xi,k)\in \widehat{\Omega}$. We decompose the frequency space $\widehat{\Omega}$ into four parts as follows.
\[
I_{1}=\left\{(\xi,k) \in \mathbb{R}\times \mathbb{N}^{+}:\xi^2< \frac{\nu^{2}}{16}(\xi^{2}+ \pi^{2}k^{2})^{3}\right\},
\]
\[
I_{2}=\left\{(\xi,k) \in \mathbb{R}\times \mathbb{N}^{+}:\frac{\nu^{2}}{16}(\xi^{2}+ \pi^{2}k^{2})^{3}\leq \xi^{2} < \frac{\nu^{2}}{4}(\xi^{2}+\pi^{2}k^{2})^{3}\right\},
\]
\[
I_{3}=\left\{(\xi,k) \in \mathbb{R}\times \mathbb{N}^{+}:\frac{\nu^{2}}{4}(\xi^{2}+\pi^{2}k^{2})^{3}\leq \xi^{2} < 4\nu^{2}(\xi^{2}+\pi^{2}k^{2})^{3}\right\},
\]
\[
I_{4}=\left\{(\xi,k) \in \mathbb{R}\times \mathbb{N}^{+}:\xi^{2}\geq 4\nu^{2}(\xi^{2}+\pi^{2}k^{2})^{3}\right\}.
\]
For $(\xi,k)\in I_{1}$, note that
\[
\lambda_{+} = -\frac{2\xi^{2}}{\nu(\xi^{2}+\pi^{2}k^{2})^{2}}\frac{1}{1+\sqrt{1-\frac{4\xi^{2}}{\nu(\xi^{2}+\pi^{2}k^{2})^{3}}}}.
\]
Thus,
\[
-\frac{2\xi^{2}}{\nu(\xi^{2}+\pi^{2}k^{2})^{2}}\leq \lambda_{+} \leq -\frac{\xi^{2}}{\nu(\xi^{2}+\pi^{2}k^{2})^{2}},\,
-\nu(\xi^{2}+\pi^{2}k^{2})\leq\lambda_{-} \leq -\frac{\nu}{2}(\xi^{2}+\pi^{2}k^{2}),
\]
\[
\frac{\nu}{2}(\xi^{2}+\pi^{2}k^{2})<\lambda_{+} -\lambda_{-} < \nu(\xi^{2}+\pi^{2}k^{2}).
\]
Since $k\geq 1$,
\begin{equation}\label{I11}
0<\widehat{\mathcal{L}_{1}(t)}\lesssim e^{-\frac{\xi^{2}}{\nu(\xi^{2}+\pi^{2}k^{2})^{2}}t},\,
0\leq\widehat{\mathcal{L}_{2}(t)}\lesssim e^{-\frac{\xi^{2}}{\nu(\xi^{2}+\pi^{2}k^{2})^{2}}t},
\end{equation}
\begin{equation}\label{I12}
\left|\widehat{\partial_{t}\mathcal{L}_{1}(t)}\right|\ls
\frac{\xi^{2}}{\nu(\xi^{2}+\pi^{2}k^{2})^{2}}e^{-\frac{\xi^{2}}{\nu(\xi^{2}+\pi^{2}k^{2})^{2}}t}
+\nu(\xi^{2}+\pi^2k^{2})e^{-\frac{\nu(\xi^{2}+\pi^{2}k^{2})}{2}t},
\end{equation}
\begin{equation}\label{I13}
\left|\widehat{\partial_{t}\mathcal{L}_{2}(t)}\right|\ls
\frac{\xi^{2}}{\nu(\xi^{2}+\pi^{2}k^{2})^{2}}e^{-\frac{\xi^{2}}{\nu(\xi^{2}+\pi^{2}k^{2})^{2}}t}
+\nu(\xi^{2}+\pi^2k^{2})e^{-\frac{\nu(\xi^{2}+\pi^{2}k^{2})}{2}t}.
\end{equation}
For $(\xi,k) \in I_{2}$,
\[
-\frac{2\xi^{2}}{\nu(\xi^{2}+\pi^{2}k^{2})^{2}}\leq \lambda_{+} \leq -\frac{\xi^{2}}{\nu(\xi^{2}+\pi^{2}k^{2})^{2}}\leq -\frac{\nu(\xi^{2}+\pi^{2}k^{2})}{16},
\]
\[
-\nu(\xi^{2}+\pi^{2}k^{2})\leq\lambda_{-} \leq -\frac{\nu}{2}(\xi^{2}+\pi^{2}k^{2}),\,
0<\lambda_{+} -\lambda_{-} <\nu(\xi^{2}+\pi^{2}k^{2}).
\]
After a straightforward calculation, we obtain that for $k\geq1$
\begin{equation}\label{I21}
0<\widehat{\mathcal{L}_{1}(t)}\lesssim e^{-\frac{\nu(\xi^2+\pi^{2}k^2)}{16}t},\,\,
\left|\widehat{\partial_{t}\mathcal{L}_{1}(t)}\right|\ls
e^{-\frac{\nu(\xi^2+\pi^2k^2)}{32}t},
\end{equation}
\begin{equation}\label{99i}
0\leq\widehat{\mathcal{L}_{2}(t)}\lesssim e^{\lambda_{+}t}t\lesssim e^{-\frac{\nu(\xi^2+\pi^{2}k^2)}{32}t},
\end{equation}
\begin{equation}\label{I22}
\left|\widehat{\partial_{t}\mathcal{L}_{2}(t)}\right|\ls \left |e^{\lambda_{+}t}\frac{\lambda_{+}-\lambda_{-}}{\lambda_{+}-\lambda_{-}} \right|+ \left|\lambda_{-}\frac{e^{\lambda_{+}t}-e^{\lambda_{-}t}}{\lambda_{+}-\lambda_{-}}\right|
\ls e^{\lambda_{+}t}+|\lambda_{-}|e^{\lambda_{+}t}t \ls
e^{-\frac{\nu(\xi^2+\pi^{2}k^2)}{32}t}.
\end{equation}
For $(\xi,k) \in I_{3}$,
\[
\sigma = i\sqrt{\frac{4\xi^{2}}{\xi^{2}+\pi^{2}k^{2}}-\nu^{2}(\xi^{2}+\pi^{2}k^{2})^{2}},\
\]
\[
\lambda_{+} = -\frac{\nu(\xi^{2}+\pi^{2}k^{2})}{2}+\frac{i}{2}\sqrt{\frac{4\xi^{2}}{\xi^{2}+\pi^{2}k^{2}}-\nu^{2}(\xi^{2}+\pi^{2}k^{2})^{2}}, \]
\[
\lambda_{-} = -\frac{\nu(\xi^{2}+\pi^{2}k^{2})}{2}-\frac{i}{2}\sqrt{\frac{4\xi^{2}}{\xi^{2}+\pi^{2}k^{2}}-\nu^{2}(\xi^{2}+\pi^{2}k^{2})^{2}}.
\]
Then
\begin{equation}\label{I31}
 \left|\widehat{\mathcal{L}_{1}(t)}\right|\lesssim e^{-\frac{\nu(\xi^{2}+\pi^{2}k^{2})}{2}t},\,
 \left|\widehat{\mathcal{L}_{2}(t)}\right|\lesssim e^{-\frac{\nu(\xi^{2}+\pi^{2}k^{2})}{2}t}\frac{\sin(\frac{|\lambda_{+}-\lambda_{-}|}{2}t)}{|\lambda_{+}-\lambda_{-}|}\lesssim e^{-\frac{\nu(\xi^{2}+\pi^{2}k^{2})}{4}t}.
\end{equation}
Moreover, since
\[
|\lambda_{+}|=|\lambda_{-}|= \frac{\xi^{2}}{\xi^{2}+\pi^{2}k^{2}} \leq 1,
\]
one has
\begin{equation}\label{I32}
\left|\widehat{\partial_{t}\mathcal{L}_{1}(t)}\right|\ls |\lambda_+|\left|e^{\lambda_+t}\right| + |\lambda_-|\left|e^{\lambda_-t}\right|\ls
e^{-\frac{\nu(\xi^{2}+\pi^{2}k^{2})}{2}t},
\end{equation}
\begin{equation}\label{I33}
\left|\widehat{\partial_{t}\mathcal{L}_{2}(t)}\right|\ls \left |e^{\lambda_{+}t}\frac{\lambda_{+}-\lambda_{-}}{\lambda_{+}-\lambda_{-}} \right|+ \left|\lambda_{-}\frac{e^{\lambda_{+}t}-e^{\lambda_{-}t}}{\lambda_{+}-\lambda_{-}}\right|
 \ls
e^{-\frac{\nu(\xi^{2}+\pi^{2}k^{2})}{4}t}.
\end{equation}
For $(\xi,k) \in I_{4}$, due to the fact that
\[
\left| \lambda_{+}-\lambda_{-} \right| = \sqrt{\frac{4\xi^{2}}{\xi^{2}+\pi^{2}k^{2}}-\nu^{2}(\xi^{2}+\pi^{2}k^{2})^{2}} \ge  \sqrt{15}\nu(\xi^{2}+\pi^2k^{2}) \gtrsim \nu,
\]
we have
\begin{equation}\label{I41}
\left|\widehat{\mathcal{L}_{1}(t)}\right|,\, \left|\widehat{\mathcal{L}_{2}(t)}\right|,\,
\left|\widehat{\partial_{t}\mathcal{L}_{1}(t)}\right|,\,
\left|\widehat{\partial_{t}\mathcal{L}_{2}}(t)\right|\ls e^{-\frac{\nu(\xi^{2}+\pi^{2}k^{2})}{2}t}.
\end{equation}
Based on estimates in (\ref{I11})-(\ref{I41}) of $\mathcal{L}_{1}$ and $\mathcal{L}_{2}$, for $0<\nu <\nu_{*} $ one can obtain (\ref{All1})-(\ref{All6}) as before, which concludes the proof of Lemma \ref{AL2}.
\hfill$\square$
\begin{Remark}
\emph{
From Lemma \ref{AL2}, we find that one derivative in the $x$-direction improves $\frac{1}{2}$-order of decay rate.
}
\end{Remark}

Using Lemma \ref{PL3}, Lemma \ref{AL1} together with Lemma \ref{AL2}, we obtain the following decay estimates for solutions to the initial-boundary value problem (\ref{Aa1})-(\ref{Aam1}).
\begin{Lemma}\label{EL2}
Assume that $\phi_{0}\in \mathfrak{D}^{8,1}, \phi_{1}\in \mathfrak{D}^{6,1}$ and $F\in L^{\infty}(0,\infty;\mathfrak{D}^{6,1})$ with
\[
\|\langle \cdot \rangle^{s} F(\cdot)\|_{L^{\infty}(W^{6,1})}:= \sup_{t>0}\langle t \rangle^{s} \|F(t)\|_{W^{6,1}} < \infty \text{ for some }s>0.
\]
Then the solution to (\ref{Aa1})-(\ref{Aam1}) satisfies for any $t>0$,
\begin{equation}\label{ALe1}
  \|\phi(t)\|_{L^{\infty}}\ls \langle t \rangle^{-\frac{1}{2}}\left(\|\phi_{0}\|_{W^{6,1}} + \|\phi_{1}\|_{W^{4,1}} + \|\langle \cdot \rangle^{s} F(\cdot)\|_{L^{\infty}(W^{4,1})}\right), \text{ if } s>1,
  \end{equation}
  \begin{equation}\label{ALe1tr}
  \|\partial_2\phi(t)\|_{L^{\infty}}\ls \langle t \rangle^{-\frac{1}{2}}\left(\|\phi_{0}\|_{W^{7,1}} + \|\phi_{1}\|_{W^{5,1}} + \|\langle \cdot \rangle^{s} F(\cdot)\|_{L^{\infty}(W^{5,1})}\right), \text{ if } s>1,
  \end{equation}
 \begin{equation}\label{ALe3}
  \|\p_{1}\phi(t)\|_{L^{\infty}}\ls \langle t \rangle^{-1}\left(\|\phi_{0}\|_{W^{8,1}} + \|\phi_{1}\|_{W^{6,1}} + \|\langle \cdot \rangle^{s} F(\cdot)\|_{L^{\infty}(W^{6,1})}\right), \text{ if } s>1,
   \end{equation}
\begin{equation}\label{ALe4}
  \|\phi(t)\|_{H^{4}}\ls \langle t \rangle^{-\frac{1}{4}}\left(\|\phi_{0}\|_{W^{8,1}} + \|\phi_{1}\|_{W^{6,1}} + \|\langle \cdot \rangle^{s} F(\cdot)\|_{L^{\infty}(W^{6,1})}\right), \text{ if } s>1,
    \end{equation}
\begin{equation}\label{ALe7}
  \|\p_{1}\phi(t)\|_{H^{2}}  \ls \langle t \rangle^{-\frac{3}{4}}\left(\|\phi_{0}\|_{W^{8,1}} + \|\phi_{1}\|_{W^{6,1}} + \|\langle \cdot \rangle^{s} F(\cdot)\|_{L^{\infty}(W^{6,1})}\right), \text{ if } s>1,
 \end{equation}
 \begin{equation}\label{ALe9+}
  \|\p_{11}\phi(t)\|_{L^{2}}\ls \langle t \rangle^{-\frac{5}{4}}\left(\|\phi_{0}\|_{W^{8,1}} + \|\phi_{1}\|_{W^{6,1}} + \|\langle \cdot \rangle^{s} F(\cdot)\|_{L^{\infty}(W^{6,1})}\right), \text{ if } s\geq\frac{5}{4},
 \end{equation}
 \begin{equation}\label{ALe9}
  \|\partial_{t}\phi(t)\|_{H^{3}}\ls \langle t\rangle^{-\frac{5}{4}}\left(\|\phi_{0}\|_{W^{8,1}} + \|\phi_{1}\|_{W^{6,1}} + \|\langle \cdot \rangle^{s} F(\cdot)\|_{L^{\infty}(W^{6,1})}\right), \text{ if } s\geq\frac{5}{4}.
\end{equation}
\end{Lemma}

\subsection{Decay rates for solutions to the linearized equation of vorticity}
Going back to (\ref{Eb1more})-(\ref{Eb2}),
we use Duhamel's principle to find
\begin{equation}\label{voint1}
w(x,y,t) = e^{\nu t\Delta}w_{0}(x,y) + \int_{0}^{t}e^{\nu (t-\tau)\Delta}\p_{1}\phi(x,y,\tau) \dtau,
\end{equation}
that is,
\begin{equation}\label{voint2}
\widehat{w}(\xi,k,t) = e^{-\nu(\xi^2+\pi^{2}k^2)t}\widehat{w_0}(\xi,k) + \int_{0}^{t}e^{-\nu(\xi^2+\pi^{2}k^2)(t-\tau)}\xi\widehat{\phi}(\xi,k,\tau)\dtau,\,(\xi,k)\in \widehat{\Omega}.
\end{equation}
\begin{Lemma}\label{EL3}
Assume that $\phi_0\in \mathfrak{D}^{8,1}$ and $\mathbf{u}_0 = (u_{10},u_{20})$ with $u_{10}\in \mathfrak{N}^{6,1}$, $u_{20}\in \mathfrak{D}^{6,1}$ so that $w_{0}\in \mathfrak{D}^{5,1}$. Then the solution $w$ to (\ref{Eb1more})-(\ref{Eb2}) satisfies
\begin {equation}\label{vor1}
 \|w(t)\|_{L^{\infty}}\ls \langle t \rangle^{-1}\left(\|w_{0}\|_{W^{5,1}} + \|\phi_{0}\|_{W^{8,1}}\right),\
  \|w(t)\|_{L^{2}}\ls \langle t \rangle^{-\frac{3}{4}}\left(\|w_{0}\|_{W^{3,1}} + \|\phi_{0}\|_{W^{6,1}}\right),
\end{equation}
\begin {equation}\label{vor2}
  \|\p_{1}w(t)\|_{L^{2}}\ls \langle t \rangle^{-\frac{5}{4}}\left(\|w_{0}\|_{W^{5,1}} + \|\phi_{0}\|_{W^{8,1}}\right),\
  \|\p_{2}w(t)\|_{L^{2}}\ls \langle t \rangle^{-\frac{3}{4}}\left(\|w_{0}\|_{W^{4,1}} + \|\phi_{0}\|_{W^{7,1}}\right),
\end {equation}
\begin{equation}\label{vor3}
\|\nabla w(t)\|_{H^{1}}\ls \langle t \rangle^{-\frac{3}{4}}\left(\|w_{0}\|_{W^{5,1}} + \|\phi_{0}\|_{W^{8,1}}\right),\
\|\p_{t}w(t)\|_{H^{2}}\ls \langle t \rangle^{-\frac{5}{4}}\left(\|w_{0}\|_{W^{5,1}} + \|\phi_{0}\|_{W^{8,1}}\right).
\end{equation}
\end{Lemma}
{\bf Proof.}
From (\ref{Eb1more})-(\ref{Eb2}), we find
\[
\phi_1(\mathbf{x}) =\partial_t \phi(\mathbf{x},0) = u_2(\mathbf{x},0)=u_{20}(\mathbf{x}) \in \mathfrak{D}^{6,1}.
\]
By Lemma \ref{EL2} and Remark \ref{Rema1},
\[
 \|\partial_1\phi(t)\|_{L^{\infty}}\ls \|\xi \widehat{\phi}(t)\|_{\widehat{L}^1} \ls \langle t \rangle^{-1}\left(\|(\xi^2 + \pi^2k^2)^3 \widehat{\phi_1}\|_{\widehat{L}^{\infty}} + \|\phi_{0}\|_{W^{8,1}}\right)
\]
\begin{equation}\label{theta1}
 \ls \langle t \rangle^{-1}\left(\|w_{0}\|_{W^{5,1}} + \|\phi_{0}\|_{W^{8,1}}\right),
\end{equation}
since $\widehat{\phi}_1(\xi,k)=\widehat{u}_{20} (\xi,k)= -\frac{i\xi}{\xi^2+\pi^2k^2}\widehat{w}_0(\xi,k)$.
By using (\ref{eLL3}), (\ref{voint2}) and (\ref{theta1}),
\[
\|w(t)\|_{L^{\infty}}\ls\|\widehat{w}(t)\|_{\widehat{L}^{1}}\ls e^{-\nu t}\|\omega_0\|_{H^2} + \int_0^t e^{-\nu (t- \tau)}\|\xi \widehat{\phi}(\tau)\|_{\widehat{L}^1}\dtau
\]
\[
\ls e^{-\nu t}\|w_0\|_{H^2} + \int_0^t e^{-\nu (t- \tau)}\langle \tau \rangle^{-1}\dtau\left(\|w_{0}\|_{W^{5,1}} + \|\phi_{0}\|_{W^{8,1}}\right)
\]
\[
\ls \langle t \rangle^{-1}\left(\|w_{0}\|_{W^{5,1}} + \|\phi_{0}\|_{W^{8,1}}\right).
\]
Other estimates in (\ref{vor1})-(\ref{vor3}) can be proved in a similar way.
\hfill$\square$

\begin{Remark}
\emph{
We observe that the decay rates of $w$ are consistent with $\partial_1 \phi$. Note that $w$ also satisfies
\[
\partial_{tt}w-\nu\Delta\partial_{t}w - \partial_{11}(-\Delta)^{-1}w = 0.
\]
}
\end{Remark}
\section{Nonlinear stability}
The local existence and uniqueness of the solution to (\ref{P4})-(\ref{P5}) can be proved by using the method of \cite{AD1,AD2}. Here we omit the details.
 \begin{Proposition}\label{tp1}
Assume that $\theta_{0}\in \mathfrak{D}^{m+1}$ and $\omega_{0} \in \mathfrak{D}^{m}$, $m\ge 2$. There exists $T^*\in (0,\infty]$ such that (\ref{P4})-(\ref{P5}) admits a unique solution
\[
(\omega, \theta) \in C([0,T^*); \mathfrak{D}^{m})\times C([0,T^*);\mathfrak{D}^{m+1}).
\]
\end{Proposition}
\begin{Remark}
\emph{
Let $(\omega, \theta)$ be a sufficiently smooth solution to (\ref{P4})-(\ref{P5}). If we assume $\theta_0\in\mathfrak{D}^{m+1}$ for some integer $m\ge 1$, then it is necessary that $\omega_0\in\mathfrak{D}^m$ and for any $t>0$,
\begin{equation}\label{P8}
\theta(t,\cdot) \in \mathfrak{D}^{m+1}, \, \omega(t,\cdot)\in\mathfrak{D}^{m},\  u_{1}(t,\cdot)\in \mathfrak{N}^{m+1},~u_{2}(t,\cdot) \in\mathfrak{D}^{m+1}.
\end{equation}
We refer to \cite{AD1,AD2,LD1} for a detailed argument.
}
\end{Remark}

Based on Proposition \ref{tp1}, global-in-time existence will follow from uniform-in-time a priori estimates. In the following we focus on decay estimates of the solutions to (\ref{P4})-(\ref{P5}).

\subsection{Integral form of solutions}
From $(\ref{pr1})$ and $(\ref{P8})$, $\theta$ and $\omega$ can be written as
\begin{equation}\label{theta2}
\theta(x,y,t) = \frac{1}{\sqrt{ \pi}}\sum_{k=1}^{+\infty}\int_{\mathbb{R}}\widehat{\theta}(\xi,k,t)e^{ i x\xi}\dki\sin k\pi y,
\end{equation}
\begin{equation}\label{vorf2}
\omega(x,y,t) = \frac{1}{\sqrt{\pi}}\sum_{k=1}^{+\infty}\int_{\mathbb{R}}\widehat{\omega}(\xi,k,t)e^{ i x\xi}\dki  \sin k\pi y.
\end{equation}
Recall that $u_{1} =\partial_{2}(-\Delta)^{-1}\omega$ and $u_{2} =-\partial_{1}(-\Delta)^{-1}\omega$, we obtain
\begin{equation}\label{vel1}
u_{1}(x,y,t) = \frac{1}{\sqrt{\pi}}\sum_{k=1}^{+\infty}\int_{\mathbb{R}}\frac{k\pi}{\xi^{2}+\pi^{2}k^{2}}\widehat{\omega}(\xi,k,t)e^{ix\xi}\dki \cos k\pi y
\end{equation}
and
\begin{equation}\label{vel2}
u_{2}(x,y,t) = \frac{1}{\sqrt{\pi}}\sum_{k=1}^{+\infty}\int_{\mathbb{R}}\frac{-i\xi }{\xi^{2}+ \pi^{2}k^{2}}\widehat{\omega}(\xi,k,t)e^{ix\xi }\dki \sin k\pi y.
\end{equation}
Taking the time derivative on $(\ref{P4})_{2}$ and using
\[
\partial_t u_2 = -\partial_1 (-\Delta)^{-1}\partial_t \omega,\,  -\Delta u_2 = \partial_t \Delta \theta + \Delta (\mathbf{u}\cdot\nabla\theta),
\]
the system $(\ref{P4})$ is transformed into
\begin{equation}\label{Eb4}
     \left\{
     \begin{aligned}
       &\partial_{t}\omega-\nu\Delta\omega = f_{1} ,\\
       &\partial_{tt}\theta-\nu\Delta\partial_{t}\theta - \partial_{1}(-\Delta)^{-1}\partial_{1}\theta = f_{2}
     \end{aligned}
     \right.
\end{equation}
with
\begin{equation}\label{Ebf4+}
f_{1} = -\mathbf{u}\cdot\nabla\omega-\p_{1}\theta,
\end{equation}
\begin{equation}\label{Ebf4++}
 f_{2} = -\partial_{t}\mathbf{u}\cdot\nabla\theta - \mathbf{u}\cdot\nabla\partial_{t}\theta - \partial_{1}(-\Delta)^{-1}(\mathbf{u}\cdot\nabla\omega)+
\nu\Delta(\mathbf{u}\cdot\nabla\theta).
\end{equation}
From (\ref{P8}), it is not difficult to find that
\[
f_1 \in \mathfrak{D}^{m-1},\, f_2 \in \mathfrak{D}^{m-2}.
\]
Accordingly,
\begin{equation}\label{Eb6}
   \omega(x,y,t)= e^{\nu t\Delta }\omega_{0}(x,y) + \int_{0}^{t}e^{\nu(t-\tau)\Delta}f_{1}(x,y,\tau)\dtau,
\end{equation}
\[
 \theta(x,y,t)=  \mathcal{L}_{1}(t)\theta_{0}(x,y) + \mathcal{L}_{2}(t)\left(\frac{\nu}{2}(-\Delta)\theta_{0}(x,y) + \theta_1(x,y)\right)
 \]
\begin{equation}\label{Eb7}
+ \int_{0}^{t}\mathcal{L}_{2}(t-\tau)f_{2}(x,y,\tau)\dtau,
\end{equation}
where
\begin{equation}\label{inithe}
\theta_1(x,y) = \p_{t}\theta(x,y,0) = -\mathbf{u}_0\cdot\nabla\theta_0 + u_{20},
\end{equation}
\begin{equation}\label{inithe2}
\mathbf{u}_0 =(u_{10}, u_{20}) = (\partial_2(-\Delta)^{-1}\omega_0, -\partial_1(-\Delta)^{-1}\omega_0).
\end{equation}
\subsection{Decay estimates of nonlinear equations}
For $m\in\mathbb{N}$ and $\epsilon>0$, which is a small parameter to be determined later, we define
\[
\mathcal{F}_{1}(t) = \sup_{0\leq \tau\leq t}\left\{\langle\tau\rangle^{-\epsilon}\left(\|\theta(\tau)\|_{H^{m+1}}+\|\omega(\tau)\|_{H^{m}}\right)\right\},
\]
\[
\mathcal{F}_{2}(t) = \sup_{0\leq \tau\leq t}\{\langle\tau\rangle\left(\|\partial_{1}\theta(\tau)\|_{L^{\infty}}+\|\mathbf{u}(\tau)\|_{W^{1,\infty}} +\|\omega(\tau)\|_{L^{\infty}}\right)
\]
\[
+\langle\tau\rangle^{\frac{1}{2}}\left(\|\partial_{2}\theta(\tau)\|_{L^{\infty}}
+\|\theta(\tau)\|_{L^{\infty}}\right)\},
\]
\begin{equation*}
\begin{split}
  \mathcal{F}_{3}(t) =& \sup_{0\leq \tau\leq t}\left\{\langle\tau\rangle^{\frac{1}{4}}\|\theta(\tau)\|_{H^{4}}
+\langle\tau\rangle^{\frac{3}{4}}\left(\|\partial_{1}\theta(\tau)\|_{H^{2}}+ \|u_{1}(\tau)\|_{H^{3}} +\|\omega(\tau)\|_{H^{2}}\right)\right. \\
  &\left.+ \langle\tau\rangle^{\frac{5}{4}}\left(\|\partial_{11}\theta(\tau)\|_{L^{2}}+\|u_{2}(\tau)\|_{H^{2}} +
   \|\p_{1}u_{1}(\tau)\|_{H^{1}} +
   \|\partial_{1}\omega(\tau)\|_{L^{2}}\right)\right\},
\end{split}
\end{equation*}
\begin{equation*}
  \mathcal{F}_4(t) =  \sup_{0\leq \tau\leq t}\left\{\langle\tau\rangle^{\frac{5}{4}}\left(\|\partial_{\tau}\theta(\tau)\|_{H^{3}}
+\|\partial_{\tau}\mathbf{u}(\tau)\|_{H^{3}}
+\|\partial_{\tau}\omega(\tau)\|_{H^{2}}\right)\right\},
\end{equation*}
\[
\mathcal{F}_{0} = \|\theta_{0}\|_{W^{8,1}} + \|\theta_{0}\|_{H^{m+1}} +\|\omega_{0}\|_{W^{5,1}} + \|\omega_{0}\|_{H^{m}},
\]
\[
\mathcal{F}(t) = \mathcal{F}_{1}(t) +\mathcal{F}_{2}(t) +\mathcal{F}_{3}(t)+\mathcal{F}_{4}(t).
\]

We show the estimates of $\mathcal{F}_j(t)$, $j=1,2,3,4$ by the following four lemmas.
\begin{Lemma}\label{EL1}
Let $m\ge 1$. Then
\[
\mathcal{F}_{1}(t) \lesssim \mathcal{F}_{0} + \left(\mathcal{F}^{3}(t) + \mathcal{F}^{4}(t)\right)^{\frac{1}{2}}.
\]
\end{Lemma}
{\bf Proof.}
The basic energy estimate for (\ref{P4}) reads as
\[
\|\nabla\theta(t)\|_{L^{2}}^{2} + \|\omega(t)\|_{L^{2}}^{2} + 2\nu\int_{0}^{t}\|\nabla\omega(\tau)\|_{L^{2}}^{2}\dtau
\]
\begin{equation}\label{Eoe2}
 \lesssim \|\nabla\theta_{0}\|_{L^{2}}^{2} + \|\omega_{0}\|_{L^{2}}^{2} + \int_{0}^{t}\|\nabla\mathbf{u}(\tau)\|_{L^{\infty}}\|\nabla\theta(\tau)\|_{L^{2}}^{2}\dtau,
\end{equation}
which follows from testing the first and second equation of (\ref{P4}) by $\omega$ and $-\Delta\theta$ respectively. Note that here we used the fact
\[
\langle \partial_1\theta,\omega \rangle + \langle -u_2, -\Delta\theta\rangle =\langle \partial_1\theta,\omega \rangle + \langle \partial_1(-\Delta)^{-1}\omega, -\Delta\theta\rangle
\]
\[
=\langle \partial_1\theta,\omega \rangle - \langle \Delta\partial_1(-\Delta)^{-1}\omega, \theta\rangle = \langle \partial_1\theta,\omega \rangle - \langle \omega, \partial_1\theta\rangle = 0.
\]
Similarly, for $m\ge 1$,
\begin{equation}\label{Ea+1}
 \langle\partial_{t}\partial^{m}\omega, \partial^{m}\omega\rangle + \langle\partial^{m}(\mathbf{u}\cdot \nabla \omega), \partial^{m}\omega\rangle - \nu\langle\Delta\partial^{m}\omega, \partial^{m}\omega\rangle = \langle\partial^{m}\partial_{1}\theta, \partial^{m}\omega\rangle
\end{equation}
and
\begin{equation}\label{Ea+2}
 \langle\partial_{t}\partial^{m}\nabla\theta, \partial^{m}\nabla\theta\rangle + \langle\partial^{m}\nabla(\mathbf{u}\cdot \nabla \theta), \partial^{m}\nabla\theta\rangle
= \langle-\partial^{m}\nabla u_{2}, \partial^{m}\nabla\theta\rangle.
\end{equation}
Adding (\ref{Ea+1}) to (\ref{Ea+2}) yields
\begin{equation}\label{Ea1}
\frac{1}{2}\frac{d}{dt}\left(\|\theta\|_{\dot{H}^{m + 1}}^{2} + \|\omega\|_{\dot{H}^{m}}^{2}\right) + \nu\|\omega\|_{\dot{H}^{m + 1}}^{2} = B_{1} + B_{2} + B_{3}
\end{equation}
with
\[
B_{1} = - \langle\partial^{m}(\mathbf{u}\cdot\nabla\omega), \partial^{m}\omega\rangle,
\]
\[
B_{2} = -\langle\partial^{m}\nabla(\mathbf{u}\cdot\nabla\theta), \partial^{m}\nabla\theta\rangle,
\]
\[
B_{3} = \langle-\partial^{m}\nabla u_{2}, \partial^{m}\nabla\theta\rangle + \langle\partial^{m}\partial_{1}\theta, \partial^{m}\omega\rangle.
\]
By the commutator estimate (\ref{eLL1}) in Lemma \ref{PL1} and Corollary \ref{corollary1},
\[
 |B_{1}|\lesssim  |\langle\partial^{m}(\mathbf{u}\cdot\nabla\omega)-\mathbf{u}\cdot\nabla\partial^{m}\omega , \partial^{m}\omega\rangle|
 \]
 \[
 \ls \|\langle\partial^{m}\text{ div }(\mathbf{u}\omega)-\mathbf{u}\cdot\nabla\partial^{m}\omega\|_{L^2}\|\partial^{m}\omega\|_{L^2}
 \]
  \begin{equation}\label{Ea2}
 \ls \left(\|\nabla\mathbf{u}\|_{L^{\infty}} + \|\omega\|_{L^{\infty}}\right)\|\mathbf{u}\|_{H^{m+1}}\|\omega\|_{H^{m}}
\ls \|\nabla\mathbf{u}\|_{L^{\infty}}\|\omega\|_{H^{m}}^{2},
\end{equation}
where we use the fact that
\[
\langle\mathbf{u}\cdot\nabla\partial^{m}\omega, \partial^{m}\omega\rangle = 0.
\]
Applying similar argument to $B_{2}$ yields
\[
|B_{2}| \lesssim \|\nabla\mathbf{u}\|_{L^{\infty}}\|\theta\|_{H^{m+1}}^{2} + \|\theta\|_{L^{\infty}}\|\mathbf{u}\|_{H^{m+2}}\|\theta\|_{H^{m+1}}
\]
\begin{equation}\label{Ea3}
\lesssim \|\nabla\mathbf{u}\|_{L^{\infty}}\|\theta\|_{H^{m+1}}^{2} + \|\theta\|_{L^{\infty}}\|\omega\|_{H^{m+1}}\|\theta\|_{H^{m+1}}.
\end{equation}
The last term $B_{3}$ vanishes once one replaces $u_{2}$ by $-\partial_{1}(-\Delta)^{-1}\omega$ and then integrates by parts twice.
\[
B_{3} = \langle \partial^{m}\nabla \partial_{1}(-\Delta)^{-1}\omega, \partial^{m}\nabla\theta\rangle + \langle\partial^{m}\partial_{1}\theta, \partial^{m}\omega\rangle = 0.
\]
By integrating (\ref{Ea1}) in time from $0$ to $t$ and summing up (\ref{Ea2}), (\ref{Ea3}),
\[
\|\theta(t)\|_{\dot{H}^{m+1}}^{2} + \|\omega(t)\|_{\dot{H}^{m}}^{2} + 2\nu\int_{0}^{t}\|\omega(\tau)\|_{\dot{H}^{m+1}}^{2}\dtau
\]
\[
 \leq  \|\theta_{0}\|_{\dot{H}^{m+1}}^{2}
+ \|\omega_{0}\|_{\dot{H}^{m}}^{2} + C\int_{0}^{t}\|\nabla\mathbf{u}(\tau)\|_{L^{\infty}}(\|\omega(\tau)\|_{H^{m}}^{2}+\|\theta(\tau)\|_{H^{m+1}}^{2})\dtau
\]
\begin{equation}\label{Eoe1}
+ \frac{C}{\nu}\int_{0}^{t}\|\theta(\tau)\|_{L^{\infty}}^{2}\|\theta(\tau)\|_{H^{m+1}}^{2}\dtau + \nu\int_{0}^{t}\|\omega(\tau)\|_{H^{m+1}}^{2}\dtau,
\end{equation}
where Young's inequality is applied.

Adding (\ref{Eoe2}) and (\ref{Eoe1}) together with Poincar$\acute{e}$ inequality,
\[
\|\theta(t)\|_{H^{m+1}}^{2} + \|\omega(t)\|_{H^{m}}^{2} + \nu\int_{0}^{t}\|\omega(\tau)\|_{H^{m+1}}^{2}\dtau \lesssim \|\theta_{0}\|_{H^{m+1}}^{2}+\|\omega_{0}\|_{H^{m}}^{2}
\]
\[
+ \int_{0}^{t}\|\nabla\mathbf{u}(\tau)\|_{L^{\infty}}\left(\|\omega(\tau)\|_{H^{m}}^{2}+\|\theta(\tau)\|_{H^{m+1}}^{2}\right)\dtau + \int_{0}^{t}\|\theta(\tau)\|_{L^{\infty}}^{2}\|\theta(\tau)\|_{H^{m+1}}^{2}\dtau.
\]
Hence,
\[
\|\theta(t)\|_{H^{m+1}}^{2} + \|\omega(t)\|_{H^{m}}^{2}
\lesssim \|\theta_{0}\|_{H^{m+1}}^{2}+\|\omega_{0}\|_{H^{m}}^{2}
\]
\[
+
\int_{0}^{t}\langle\tau\rangle^{-1}\langle\tau\rangle^{2\epsilon}\langle\tau\rangle\|\nabla\mathbf{u}(\tau)\|_{L^{\infty}}\dtau \sup_{0\leq \tau\leq t}\left\{\langle\tau\rangle^{-2\epsilon}\|\theta(\tau)\|_{H^{m+1}}^{2}+\langle\tau\rangle^{-2\epsilon}\|\omega(\tau)\|_{H^{m}}^{2}\right\}
\]
\[
+\int_{0}^{t}\langle\tau\rangle^{-1}\langle\tau\rangle\|\theta(\tau)\|_{L^{\infty}}^{2}\langle\tau\rangle^{2\epsilon}\dtau \sup_{0\leq \tau\leq t}\left\{\langle\tau\rangle^{-2\epsilon}\|\theta(\tau)\|_{H^{m+1}}^{2}+\langle\tau\rangle^{-2\epsilon}\|\omega(\tau)\|_{H^{m}}^{2}\right\}
\]
\[
\lesssim \mathcal{F}_{0}^{2} + \mathcal{F}_{1}^{2}(t)\mathcal{F}_{2}(t)\int_{0}^{t}\langle\tau\rangle^{-1}\langle\tau\rangle^{2\epsilon}\dtau +
\mathcal{F}_{1}^{2}(t)\mathcal{F}_{2}^{2}(t)\int_{0}^{t}\langle\tau\rangle^{-1}\langle\tau\rangle^{2\epsilon}\dtau
\]
\[
\lesssim \mathcal{F}_{0}^{2} + \langle t\rangle^{2\epsilon}\mathcal{F}^{3}(t)+ \langle t\rangle^{2\epsilon}\mathcal{F}^{4}(t).
\]
Finally, we conclude that
\[
\mathcal{F}_{1}(t) \lesssim \mathcal{F}_{0} + \left(\mathcal{F}^{3}(t) + \mathcal{F}^{4}(t)\right)^{\frac{1}{2}}.
\]
\hfill$\square$
\begin{Lemma}\label{EL4}
Let $m> 20$. Then for $0<\delta <\frac{2}{5}$,
\[
\mathcal{F}_{2}(t) \lesssim \mathcal{F}_{0} + \mathcal{F}_{0}^{2} + \mathcal{F}^{2}(t) + \mathcal{F}^{\frac{5}{2}}(t) + \mathcal{F}^{2 + \delta}(t).
\]
\end{Lemma}
{\bf Proof.}
From (\ref{Eb7}) and Lemma \ref{AL1}-\ref{AL2}, we obtain that
\[
\|\partial_{1}\theta(t)\|_{L^{\infty}} \lesssim \|\widehat{\partial_{1}\theta}(t)\|_{\widehat{L}^{1}}
\]
\[
\ls
\langle t\rangle^{-1}\|\omega_{0}\|_{W^{5,1}} + \langle t\rangle^{-1}\|\theta_{0}\|_{W^{8,1}} +\langle t\rangle^{-1}\|\theta_{0}\|_{H^{m}}\|\omega_{0}\|_{H^{m}}
+ J_{1} + J_{2} + J_{3} + J_{4},
\]
where
\[
J_{1}=\int_{0}^{t}\langle t-\tau \rangle^{-1}\|\left(\partial_{\tau}\mathbf{u}\cdot\nabla\theta\right)(\tau)\|_{W^{6,1}}\dtau,\
J_{2}=\int_{0}^{t}\langle t-\tau\rangle^{-1}\|\left(\mathbf{u}\cdot\nabla\partial_{\tau}\theta\right)(\tau)\|_{W^{6,1}}\dtau,
\]
\[
J_{3}=\int_{0}^{t}\langle t-\tau\rangle^{-1}\|\partial_{1}(\mathbf{u}\cdot\nabla\omega)(\tau)\|_{W^{4,1}}\dtau,\
J_{4}=\int_{0}^{t}\langle t-\tau\rangle^{-1}\|\left(\mathbf{u}\cdot\nabla\theta\right)(\tau)\|_{W^{8,1}}\dtau.
\]

For $J_{1}$, according to the interpolation inequality (\ref{BLT1}) together with (\ref{BL13}) in Lemma \ref{PL1},
\[
\|\partial_{t}\mathbf{u}\cdot\nabla\theta\|_{W^{6,1}}
\ls \|\partial_{t}\mathbf{u}\|_{L^{2}}\|\nabla\theta\|_{H^{6}} + \|\partial_{t}\mathbf{u}\|_{H^{6}}\|\nabla\theta\|_{L^{2}}
\]
\[
\ls \|\partial_{t}\mathbf{u}\|_{L^{2}}\|\nabla\theta\|_{H^{12}}^{\frac{1}{2}}\|\nabla\theta\|_{L^{2}}^{\frac{1}{2}} + \|\partial_{t}\mathbf{u}\|_{H^{\frac{6}{\delta}}}^{\delta}\|\partial_{t}\mathbf{u}\|_{L^{2}}^{1-\delta}\|\nabla\theta\|_{L^{2}},
\]
where $0<\delta  < \frac{2}{5}$. Replacing $\|\partial_{t}\mathbf{u}\|_{H^{\frac{6}{\delta}}}$ by $\|\partial_{t}\omega\|_{H^{\frac{6}{\delta}-1}}$ and using $(\ref{P4})_{1}$, one infers that
\[
J_{1}\ls \int_{0}^{t}\langle t-\tau\rangle^{-1}\langle \tau\rangle^{-\frac{5}{4}}\langle \tau\rangle^{\frac{5}{4}}\|\partial_{\tau}\mathbf{u}\|_{L^{2}}\langle t\rangle^{-\frac{\epsilon}{2}}\langle \tau \rangle^{\frac{\epsilon}{2}}\|\nabla\theta\|_{H^{12}}^{\frac{1}{2}}
\langle\tau\rangle^{-\frac{1}{8}}\langle\tau\rangle^{\frac{1}{8}}
\|\nabla\theta\|_{L^{2}}^{\frac{1}{2}}\dtau
\]
\[
+ \int_{0}^{t}\langle t-\tau\rangle^{-1}
\langle \tau \rangle^{-\delta\epsilon}\langle \tau \rangle^{\delta\epsilon}\|\Delta\omega+\partial_{1}\theta\|_{H^{\frac{6}{\delta}-1}}^{\delta}
\langle \tau\rangle^{-\frac{5}{4}(1-\delta)}\langle \tau\rangle^{\frac{5}{4}(1-\delta)}\|\p_{\tau}\mathbf{u}\|_{L^{2}}^{1-\delta}
\langle\tau\rangle^{-\frac{1}{4}}\langle\tau\rangle^{\frac{1}{4}}\|\nabla\theta\|_{L^{2}}\dtau
\]
\[
+ \int_{0}^{t}\langle t-\tau\rangle^{-1}
\langle \tau \rangle^{-2\delta\epsilon}\langle \tau \rangle^{2\delta\epsilon}\|\mathbf{u}\cdot\nabla\omega\|_{H^{\frac{\delta}{6}-1}}^{\delta}
\langle \tau\rangle^{-\frac{5}{4}(1-\delta)}
\langle\tau\rangle^{\frac{5}{4}(1-\delta)}\|\p_{\tau}\mathbf{u}\|_{L^{2}}^{1-\delta}
\langle\tau\rangle^{-\frac{1}{4}}\langle\tau\rangle^{\frac{1}{4}}\|\nabla\theta\|_{L^{2}}\dtau
\]
\[
\ls
\int_{0}^{t}\langle t-\tau\rangle^{-1}\langle
\tau\rangle^{-\frac{11}{8}+\frac{\epsilon}{2}}\dtau\mathcal{F}_{4}(t)\mathcal{F}_{1}^{\frac{1}{2}}(t)\mathcal{F}_{3}^{\frac{1}{2}}(t)
+\int_{0}^{t}\langle t-\tau\rangle^{-1} \langle \tau\rangle^{-\frac{3}{2}+\frac{5\delta}{4}+\delta\epsilon}
\dtau\mathcal{F}_{4}^{1-\delta}(t)\mathcal{F}_{1}^{\delta}(t)\mathcal{F}_{3}(t)
\]
\[
+\int_{0}^{t}\langle t-\tau\rangle^{-1}
\langle\tau\rangle^{-\frac{3}{2}+\frac{5\delta}{4}+2\delta\epsilon} \dtau\mathcal{F}_{4}^{1-\delta}(t)\mathcal{F}_{1}^{2\delta}(t)\mathcal{F}_{3}(t)
\]
\begin{equation}\label{JJ1}
 \ls \langle t\rangle^{-1}\left(\mathcal{F}^{2}(t) + \mathcal{F}^{2+\delta}(t)\right).
\end{equation}
Here $\epsilon >0$ is such that $\frac{5\delta}{4}+2\delta\epsilon < \frac{1}{2}$ according to Lemma \ref{EL2}.

Similarly for $J_{2}$, we have
\[
\|\mathbf{u}\cdot\nabla\partial_{t}\theta\|_{W^{6,1}}
\ls \|\mathbf{u}\|_{L^{2}}\|\partial_{t}\nabla\theta\|_{H^{6}} + \|\mathbf{u}\|_{H^{6}}\|\partial_{t}\nabla\theta\|_{L^{2}}
\]
\[
\ls \|\mathbf{u}\|_{L^{2}}\|\partial_{t}\nabla\theta\|_{H^{12}}^{\frac{1}{2}}\|\partial_{t}\nabla\theta\|_{L^{2}}^{\frac{1}{2}} + \|\mathbf{u}\|_{H^{12}}^{\frac{1}{2}}\|\mathbf{u}\|_{L^{2}}^{\frac{1}{2}}\|\partial_{t}\nabla\theta\|_{L^{2}}.
\]
By $(\ref{P4})_{2}$ and Lemma \ref{EL2},
\[
J_{2}\ls \int_{0}^{t}\langle t-\tau\rangle^{-1}\langle \tau \rangle^{-\frac{3}{4}}\langle \tau\rangle^{\frac{3}{4}}
\|\mathbf{u}\|_{L^{2}}\langle \tau \rangle^{-\epsilon}\langle \tau \rangle^{\epsilon}\|\mathbf{u}\cdot \nabla \theta\|_{H^{13}}^{\frac{1}{2}}\langle \tau\rangle^{-\frac{5}{8}}\langle \tau\rangle^{\frac{5}{8}}\|\partial_{\tau}\nabla\theta\|_{L^{2}}^{\frac{1}{2}}\dtau
\]
\[
+ \int_{0}^{t}\langle t-\tau\rangle^{-1}\langle \tau \rangle^{-\frac{3}{4}}\langle \tau\rangle^{\frac{3}{4}}
\|\mathbf{u}\|_{L^{2}}\langle \tau \rangle^{-\frac{\epsilon}{2}}\langle \tau \rangle^{\frac{\epsilon}{2}}\|u_{2}\|_{H^{13}}^{\frac{1}{2}}\langle \tau\rangle^{-\frac{5}{8}}\langle \tau\rangle^{\frac{5}{8}}\|\partial_{\tau}\nabla\theta\|_{L^{2}}^{\frac{1}{2}}\dtau
\]
\[
+\int_{0}^{t}\langle t-\tau\rangle^{-1}\langle \tau \rangle^{-\frac{3}{8}}\langle \tau\rangle^{\frac{3}{8}}
\|\mathbf{u}\|_{L^{2}}^{\frac{1}{2}}\langle \tau \rangle^{\frac{\epsilon}{2}}\langle \tau \rangle^{-\frac{\epsilon}{2}}\|\mathbf{u}\|_{H^{12}}^{\frac{1}{2}}\langle \tau\rangle^{-\frac{5}{4}}\langle \tau\rangle^{\frac{5}{4}}\|\partial_{\tau}\nabla\theta\|_{L^{2}}\dtau
\]
\[
\ls
\int_{0}^{t}\langle t-\tau\rangle^{-1}\langle \tau\rangle^{-\frac{11}{8}+\epsilon}\dtau \mathcal{F}_{3}(t)\mathcal{F}_{1}(t)\mathcal{F}_{4}^{\frac{1}{2}}(t)
+ \int_{0}^{t}\langle t-\tau\rangle^{-1}\langle \tau\rangle^{-\frac{11}{8}+\frac{\epsilon}{2}}\dtau
\mathcal{F}_{3}(t)\mathcal{F}_{1}^{\frac{1}{2}}(t)\mathcal{F}_{4}^{\frac{1}{2}}(t)
\]
\[
+ \int_{0}^{t}\langle t-\tau\rangle^{-1}\langle \tau\rangle^{-\frac{13}{8}+\frac{\epsilon}{2}}\dtau
\mathcal{F}_{3}^{\frac{1}{2}}(t)\mathcal{F}_{1}^{\frac{1}{2}}(t)\mathcal{F}_{4}(t)
\]
\begin{equation}\label{JJ2}
\ls \langle t\rangle^{-1}\left(\mathcal{F}^{2}(t) + \mathcal{F}^{\frac{5}{2}}(t)\right).
\end{equation}

For $J_{3}$, it holds that
\[
\|\partial_{1}(\mathbf{u}\cdot\nabla\omega)\|_{W^{4,1}} \ls \|\mathbf{u}\cdot\nabla\omega\|_{W^{5,1}}
\ls \|u_{1}\partial_{1}\omega\|_{W^{5,1}}
+ \|u_{2}\partial_{2}\omega\|_{W^{5,1}}
\]
\[
\ls \|u_{1}\|_{L^{2}}\|\partial_{1}\omega\|_{L^{2}}^{\frac{1}{2}}\|\partial_{1}\omega\|_{H^{10}}^{\frac{1}{2}} + \|u_{1}\|_{L^{2}}^{\frac{1}{2}}\|u_{1}\|_{H^{10}}^{\frac{1}{2}}\|\partial_{1}\omega\|_{L^{2}}
\]
\[
+ \|u_{2}\|_{L^{2}}\|\partial_{2}\omega\|_{L^{2}}^{\frac{1}{2}}\|\partial_{2}\omega\|_{H^{10}}^{\frac{1}{2}} + \|u_{2}\|_{L^{2}}^{\frac{1}{2}}\|u_{2}\|_{H^{10}}^{\frac{1}{2}}\|\partial_{2}\omega\|_{L^{2}}.
\]
In a similar way,
\begin{equation}\label{JJ3}
J_{3}
\ls
\int_{0}^{t}\langle t-\tau\rangle^{-1}\left(\langle \tau \rangle^{-\frac{11}{8}+\frac{\epsilon}{2}} +
\langle \tau \rangle^{-\frac{13}{8}+\frac{\epsilon}{2}}\right)\dtau\mathcal{F}^{2}(t)
\ls \langle t \rangle^{-1}\mathcal{F}^{2}(t).
\end{equation}

For the last term $J_{4}$, using the interpolation inequality (\ref{BLT1}) and (\ref{BL13}) in Lemma \ref{PL1}, we obtain that for all $\delta \in (0,1)$,
\[
\|\mathbf{u}\cdot\nabla\theta\|_{W^{8,1}}
\ls \|u_{1}\|_{L^{2}}\|\partial_{1}\theta\|_{H^{8}} + \|u_{1}\|_{H^{8}}\|\partial_{1}\theta\|_{L^{2}}
\]
\[
+ \|u_{2}\|_{L^{2}}\|\partial_{2}\theta\|_{H^{8}} + \|u_{2}\|_{H^{8}}\|\partial_{2}\theta\|_{L^{2}}
\]
\[
\ls \|u_{1}\|_{L^{2}}\|\partial_{1}\theta\|_{H^{16}}^{\frac{1}{2}}\|\partial_{1}\theta\|_{L^{2}}^{\frac{1}{2}} + \|u_{1}\|_{H^{16}}^{\frac{1}{2}}\|u_{1}\|_{L^{2}}^{\frac{1}{2}}\|\partial_{1}\theta\|_{L^{2}}
\]
\[
+ \|u_{2}\|_{L^{2}}\|\partial_{2}\theta\|_{H^{16}}^{\frac{1}{2}}\|\partial_{2}\theta\|_{L^{2}}^{\frac{1}{2}} + \|u_{2}\|_{H^{\frac{8}{\delta}}}^{\delta}\|u_{2}\|_{L^{2}}^{1-\delta}\|\partial_{2}\theta\|_{L^{2}}.
\]
By choosing $0<\delta < \frac{2}{5}$,
\begin{equation}\label{JJ4}
 J_{4} \ls
\int_{0}^{t}\langle t-\tau\rangle^{-1}\left(\langle \tau\rangle^{-\frac{9}{8}+\frac{\epsilon}{2}}+\langle \tau \rangle^{-\frac{11}{8}+\frac{\epsilon}{2}}+\langle \tau\rangle^{-\frac{3}{2}+\frac{5\delta}{4}+\delta\epsilon}\right)\dtau \mathcal{F}^{2}(t)
\ls \langle t\rangle^{-1}\mathcal{F}^{2}(t).
\end{equation}

Summing up (\ref{JJ1})-(\ref{JJ4}) gives
\begin{equation}\label{Eaa}
\|\partial_{1}\theta(t)\|_{L^{\infty}} \lesssim \|\widehat{\p_{1}\theta}(t)\|_{\widehat{L}^{1}} \lesssim
\langle t\rangle^{-1}\left(\mathcal{F}_{0} +\mathcal{F}_{0}^{2} +\mathcal{F}^{2}(t)+\mathcal{F}^{2 + \delta}(t) +
\mathcal{F}^{\frac{5}{2}}(t)\right).
\end{equation}
Note that
\[
\|\theta(t)\|_{L^{\infty}} +\|\partial_2\theta(t)\|_{L^{\infty}}\lesssim \|\widehat{\theta}(t)\|_{\widehat{L}^{1}}
+
\|\widehat{\partial_2\theta}(t)\|_{\widehat{L}^{1}}
\]
\[
 \lesssim
\langle t\rangle^{-\frac{1}{2}}\|\omega_{0}\|_{W^{4,1}} + \langle t\rangle^{-\frac{1}{2}}\|\theta_{0}\|_{W^{7,1}} +
\langle t\rangle^{-\frac{1}{2}}\|\theta_{0}\|_{H^{m}}
\|\omega_{0}\|_{H^{m}}
\]
\[
+ \int_{0}^{t}\langle t-\tau \rangle^{-\frac{1}{2}}\left(\|\left(\partial_{\tau}\mathbf{u}\cdot\nabla\theta\right)(\tau)\|_{W^{5,1}}
+ \|\left(\mathbf{u}\cdot\nabla\partial_{\tau}\theta\right)(\tau)\|_{W^{5,1}}\right)\dtau
\]
\[
+\int_{0}^{t}\langle t-\tau \rangle^{-\frac{1}{2}}\left(
 \|\partial_{1}(\mathbf{u}\cdot\nabla\omega)(\tau)\|_{W^{3,1}}
+ \|\left(\mathbf{u}\cdot\nabla\theta\right)(\tau)\|_{W^{7,1}}\right)\dtau.
\]
Using the method similar to estimate $\|\partial_{1}\theta\|_{L^{\infty}}$,
\[
\|\theta(t)\|_{L^{\infty}}+\|\partial_{1}\theta(t)\|_{L^{\infty}} \lesssim \|\widehat{\theta}(t)\|_{\widehat{L}^{1}}+\|\widehat{\partial_{1}\theta}(t)\|_{\widehat{L}^{1}} \]
\[
\ls
\langle t\rangle^{-\frac{1}{2}}\left(\mathcal{F}_{0} +\mathcal{F}_{0}^{2} +\mathcal{F}^{2}(t)+\mathcal{F}^{2 + \delta}(t) +
\mathcal{F}^{\frac{5}{2}}(t)\right).
\]

Now we turn to estimate $\omega$. By (\ref{Eb6}) and (\ref{Eaa}),
\[
\|\omega(t)\|_{L^{\infty}} \lesssim \|\widehat{\omega}(t)\|_{\widehat{L}^{1}}
\]
\[
\lesssim \left\|e^{-\nu\left(\xi^{2}+ \pi^2 k^{2}\right) t}\widehat{\omega}_{0}\right\|_{\widehat{L}^{1}} +
\int_{0}^{t}\left\|e^{-\nu\left(\xi^{2}+ \pi^2k^{2}\right)(t-\tau)}\widehat{\mathbf{u}\cdot\nabla\omega}(\tau)\right\|_{\widehat{L}^{1}} \dtau
\]
\[
+ \int_{0}^{t}\left\|e^{-\nu\left(\xi^{2}+ \pi^2k^{2}\right)(t-\tau)}\widehat{\p_{1}\theta}(\tau)\right\|_{\widehat{L}^{1}}\dtau
\]
\[
\lesssim e^{-\nu t}\|\widehat{\omega}_{0}\|_{\widehat{L}^{1}} +
\int_{0}^{t}e^{-\nu(t-\tau)}\|\widehat{\mathbf{u}\cdot\nabla\omega}(\tau)\|_{\widehat{L}^{1}}\dtau
+
\int_{0}^{t}e^{-\nu(t-\tau)}\langle\tau\rangle^{-1}\langle\tau\rangle\|\widehat{\p_{1}\theta}(\tau)\|_{\widehat{L}^{1}}\dtau
\]
\[
\lesssim e^{-\nu t}\|\omega_{0}\|_{H^{2}} +
\int_{0}^{t}e^{-\nu(t-\tau)}\|\left(\mathbf{u}\cdot\nabla\omega\right)(\tau)\|_{H^{2}}\dtau
\]
\begin{equation}\label{JJJ1}
+\langle t\rangle^{-1}\left(\mathcal{F}_{0} +
\mathcal{F}_{0}^{2} + \mathcal{F}^{2}(t)
 +\mathcal{F}^{2+\delta}(t) +
\mathcal{F}^{\frac{5}{2}}(t)\right),
\end{equation}
where (\ref{eLL3}) and the fact that $k\geq1$ are used.
It remains to estimate
\[
\int_{0}^{t}e^{-\nu(t-\tau)}\|\left(\mathbf{u}\cdot\nabla\omega\right)(\tau)\|_{H^{2}}\dtau.
\]
Note that
\[
\|\mathbf{u}\cdot\nabla\omega\|_{H^{2}} \ls \|\mathbf{u}\cdot\nabla\omega\|_{W^{3,1}}
\ls \|u_{1}\partial_{1}\omega\|_{W^{3,1}} + \|u_{2}\omega\|_{W^{4,1}}
\]
\[
\ls \|u_{1}\|_{L^{2}}\|\partial_{1}\omega\|_{L^{2}}^{\frac{1}{2}}\|\partial_{1}\omega\|_{H^{6}}^{\frac{1}{2}} + \|u_{1}\|_{L^{2}}^{\frac{1}{2}}\|u_{1}\|_{H^{6}}^{\frac{1}{2}}\|\partial_{1}\omega\|_{L^{2}}
\]
\[
+ \|u_{2}\|_{L^{2}}\|\omega\|_{L^{2}}^{\frac{1}{2}}\|\omega\|_{H^{8}}^{\frac{1}{2}} + \|u_{2}\|_{L^{2}}^{\frac{1}{2}}\|u_{2}\|_{H^{8}}^{\frac{1}{2}}\|\omega\|_{L^{2}}.
\]
Then
\[
\int_{0}^{t}e^{-\nu(t-\tau)}\|\left(\mathbf{u}\cdot\nabla\omega\right)(\tau)\|_{H^{2}}\dtau
\]
\begin{equation}\label{JJJ2}
\ls
\int_{0}^{t}e^{-\nu(t-\tau)}\left(\langle \tau\rangle^{-\frac{11}{8}+\frac{\epsilon}{2}} + \langle \tau\rangle^{-\frac{13}{8}+\frac{\epsilon}{2}}\right)\dtau\mathcal{F}^{2}(t)
\ls \langle t\rangle^{-\frac{11}{8}+\frac{\epsilon}{2}}\mathcal{F}^{2}(t).
\end{equation}
From (\ref{JJJ1}) and (\ref{JJJ2}), we find
\begin{equation}\label{ds}
\|\omega(t)\|_{L^{\infty}} \lesssim \|\widehat{\omega}(t)\|_{\widehat{L}^{1}}
\lesssim
\langle t\rangle^{-1}\left(\mathcal{F}_{0} +
\mathcal{F}_{0}^{2} + \mathcal{F}^{2}(t) +\mathcal{F}^{2+\delta}(t) +
\mathcal{F}^{\frac{5}{2}}(t)\right).
\end{equation}
According to (\ref{vel1})-(\ref{vel2}) and (\ref{ds}) together with the fact that $k\geq1$,
\[
\|\mathbf{u}(t)\|_{L^{\infty}}
\ls \left\|\frac{1}{(\xi^{2}+\pi^2k^{2})^{\frac{1}{2}}}\widehat{\omega}(t)\right\|_{\widehat{L}^{1}}
\lesssim \|\widehat{\omega}(t)\|_{\widehat{L}^{1}}
\]
\[
\lesssim
\langle t\rangle^{-1}\left(\mathcal{F}_{0} +
\mathcal{F}_{0}^{2} + \mathcal{F}^{2}(t) +\mathcal{F}^{2+\delta}(t) +
\mathcal{F}^{\frac{5}{2}}(t)\right),
\]
\[
\|\nabla\mathbf{u}(t)\|_{L^{\infty}}\lesssim \|\widehat{\omega}(t)\|_{\widehat{L}^{1}}
\lesssim
\langle t\rangle^{-1}\left(\mathcal{F}_{0} +
\mathcal{F}_{0}^{2} + \mathcal{F}^{2}(t) +\mathcal{F}^{2+\delta}(t) +
\mathcal{F}^{\frac{5}{2}}(t)\right).
\]
Thus, we finish the proof of Lemma \ref{EL4}.
\hfill$\square$
\begin{Lemma}\label{EL5}
Let $m>32$. Then for $0<\delta <\frac{1}{5}$,
\[
\mathcal{F}_{3}(t) \lesssim \mathcal{F}_{0} + \mathcal{F}_{0}^{2} + \mathcal{F}^{2}(t) + \mathcal{F}^{\frac{5}{2}}(t)
+\mathcal{F}^{2+\delta}(t).
\]
\end{Lemma}
{\bf Proof.}
One deduces from Lemma \ref{AL1}-\ref{AL2} that
\[
\|\theta(t)\|_{H^{4}} = \|\widehat{\theta}(t)\|_{\widehat{H}^{4}} = \|\left(1+ (\xi^2 + \pi^2k^2)\right)^{2}\widehat{\theta}(t)\|_{\widehat{L}^{2}}
\]
\[
\ls
\langle t\rangle^{-\frac{1}{4}}\|\omega_{0}\|_{W^{5,1}} + \langle t\rangle^{-\frac{1}{4}}\|\theta_{0}\|_{W^{8,1}} +
\langle t\rangle^{-\frac{1}{4}}\|\theta_{0}\|_{H^{m}}
\|\omega_{0}\|_{H^{m}}
\]
\[
+ \int_{0}^{t}\langle t-\tau \rangle^{-\frac{1}{4}}\left(\|\left(\partial_{\tau}\mathbf{u}\cdot\nabla\theta\right)(\tau)\|_{W^{6,1}}
+ \|\left(\mathbf{u}\cdot\nabla\partial_{\tau}\theta\right)(\tau)\|_{W^{6,1}}\right)\dtau
\]
\[
+\int_{0}^{t}\langle t-\tau \rangle^{-\frac{1}{4}}\left( \|\partial_{1}(\mathbf{u}\cdot\nabla\omega)(\tau)\|_{W^{4,1}}
+ \|\left(\mathbf{u}\cdot\nabla\theta\right)(\tau)\|_{W^{8,1}}\right)\dtau,
\]
\[
\|\partial_{1}\theta(t)\|_{H^{2}} = \|\widehat{\partial_{1}\theta}(t)\|_{\widehat{H}^{2}} \ls
\langle t\rangle^{-\frac{3}{4}}\|\omega_{0}\|_{W^{5,1}} + \langle t\rangle^{-\frac{3}{4}}\|\theta_{0}\|_{W^{8,1}} +
\langle t\rangle^{-\frac{3}{4}}\|\theta_{0}\|_{H^{m}}
\|\omega_{0}\|_{H^{m}}
\]
\[
+ \int_{0}^{t}\langle t-\tau\rangle^{-\frac{3}{4}}\left(\|\left(\partial_{\tau}\mathbf{u}\cdot\nabla\theta\right)(\tau)\|_{W^{6,1}}
+ \|\left(\mathbf{u}\cdot\nabla\partial_{\tau}\theta\right)(\tau)\|_{W^{6,1}}\right)\dtau
\]
\[
+ \int_{0}^{t}\langle t-\tau\rangle^{-\frac{3}{4}}\left(\|\partial_{1}(\mathbf{u}\cdot\nabla\omega)(\tau)\|_{W^{4,1}}
+ \|\left(\mathbf{u}\cdot\nabla\theta\right)(\tau)\|_{W^{8,1}}\right)\dtau,
\]
Similar to the  $L^\infty$-estimate of $\partial_1 \theta$ in Lemma \ref{EL4}, we obtain for $0<\delta <\frac{2}{5}$,
\begin{equation}\label{xb2}
\|\theta(t)\|_{H^{4}} = \|\widehat{\theta}(t)\|_{\widehat{H}^{4}}\ls
\langle t\rangle^{-\frac{1}{4}}\left(\mathcal{F}_{0} +
\mathcal{F}_{0}^{2}+\mathcal{F}^{2}(t) +\mathcal{F}^{\frac{5}{2}}(t)+ \mathcal{F}^{2+\delta}(t)\right),
\end{equation}
\begin{equation}\label{xb1}
 \|\partial_{1}\theta(t)\|_{H^{2}}= \|\widehat{\partial_{1}\theta}(t)\|_{\widehat{H}^{2}}  \ls
\langle t\rangle^{-\frac{3}{4}}\left(\mathcal{F}_{0} +
\mathcal{F}_{0}^{2} +
\mathcal{F}^{2}(t) + \mathcal{F}^{\frac{5}{2}}(t)
+\mathcal{F}^{2+\delta}(t)\right).
\end{equation}
Moreover,
\[
\|\partial_{11}\theta(t)\|_{L^{2}} =
\|\widehat{\partial_{11}\theta}(t)\|_{\widehat{L}^{2}} \lesssim
\langle t\rangle^{-\frac{5}{4}}\|\omega_{0}\|_{W^{5,1}} + \langle t\rangle^{-\frac{5}{4}}\|\theta_{0}\|_{W^{8,1}} +
\langle t\rangle^{-\frac{5}{4}}\|\theta_{0}\|_{H^{m}}
\|\omega_{0}\|_{H^{m}}
\]
\[
+ O_{1} + O_{2} + O_{3} + O_{4},
\]
where
\[
O_{1} = \int_{0}^{t}\langle t-\tau \rangle^{-\frac{5}{4}}\|\left(\partial_{\tau}\mathbf{u}\cdot\nabla\theta\right)(\tau)\|_{W^{6,1}}\dtau,\,
O_{2} = \int_{0}^{t}\langle t-\tau\rangle^{-\frac{5}{4}}\|\left(\mathbf{u}\cdot\nabla\partial_{\tau}\theta\right)(\tau)\|_{W^{6,1}}\dtau,
\]
\[
O_{3} = \int_{0}^{t}\langle t-\tau\rangle^{-\frac{5}{4}}\|\partial_{1}(\mathbf{u}\cdot\nabla\omega)(\tau)\|_{W^{4,1}}\dtau,\,
O_{4} = \int_{0}^{t}\langle t-\tau\rangle^{-\frac{5}{4}}\|\left(\mathbf{u}\cdot\nabla\theta\right)(\tau)\|_{W^{8,1}}\dtau.
\]
Note that
\[
 \|\partial_{t}\mathbf{u}\cdot\nabla\theta\|_{W^{6,1}}
\ls \|\partial_{t}\mathbf{u}\|_{L^{2}}\|\nabla\theta\|_{H^{6}} + \|\partial_{t}\mathbf{u}\|_{H^{6}}\|\nabla\theta\|_{L^{2}}
\]
\[
\ls \|\partial_{t}\mathbf{u}\|_{L^{2}}\|\nabla\theta\|_{H^{12}}^{\frac{1}{2}}\|\nabla\theta\|_{L^{2}}^{\frac{1}{2}} + \|\partial_{t}\mathbf{u}\|_{H^{\frac{6}{\delta}}}^{\delta}\|\partial_{t}\mathbf{u}\|_{L^{2}}^{1-\delta}\|\nabla\theta\|_{L^{2}}.
\]
Thus
\[
O_{1}
\ls
\int_{0}^{t}\langle t-\tau\rangle^{-\frac{5}{4}}\langle
\tau\rangle^{-\frac{11}{8}+\frac{\epsilon}{2}}\dtau\mathcal{F}_{4}(t)\mathcal{F}_{1}^{\frac{1}{2}}(t)\mathcal{F}_{3}^{\frac{1}{2}}(t)
\]
\[
+\int_{0}^{t}\langle t-\tau\rangle^{-\frac{5}{4}} \langle \tau\rangle^{-\frac{3}{2}+\frac{5\delta}{4}+\delta\epsilon}
\dtau\mathcal{F}_{4}^{1-\delta}(t)\mathcal{F}_{1}^{\delta}(t)\mathcal{F}_{3}(t)
\]
\begin{equation}\label{DJD1}
+\int_{0}^{t}\langle t-\tau\rangle^{-\frac{5}{4}}
\langle\tau\rangle^{-\frac{3}{2}+\frac{5\delta}{4}+2\delta\epsilon} \dtau\mathcal{F}_{4}^{1-\delta}(t)\mathcal{F}_{1}^{2\delta}(t)\mathcal{F}_{3}(t)
\ls \langle t\rangle^{-\frac{5}{4}}\left(\mathcal{F}^{2}(t) + \mathcal{F}^{2+\delta}(t)\right),
\end{equation}
where $0<\delta  < \frac{1}{5}$ and $\frac{5\delta}{4}+2\delta\epsilon < \frac{1}{4}$.
Similarly,
\[
O_{2}
\ls
\int_{0}^{t}\langle t-\tau\rangle^{-\frac{5}{4}}\left(\langle \tau\rangle^{-\frac{11}{8}+\frac{\epsilon}{2}} +\langle \tau\rangle^{-\frac{13}{8}+\frac{\epsilon}{2}}\right)\dtau
\mathcal{F}_{3}(t)\mathcal{F}_{1}(t)\mathcal{F}_{4}^{\frac{1}{2}}(t)
\]
\begin{equation}\label{DJD2}
+ \int_{0}^{t}\langle t-\tau\rangle^{-\frac{5}{4}}\langle \tau\rangle^{-\frac{11}{8}+\epsilon}\dtau\mathcal{F}_{4}(t)\mathcal{F}_{1}^{\frac{1}{2}}(t)\mathcal{F}_{3}^{\frac{1}{2}}(t)
\ls \langle t\rangle^{-\frac{5}{4}}\left(\mathcal{F}^{2}(t) + \mathcal{F}^{\frac{5}{2}}(t)\right),
\end{equation}
and
\begin{equation}\label{DJD3}
O_{3}
\ls
\int_{0}^{t}\langle t-\tau\rangle^{-\frac{5}{4}}\left(\langle \tau \rangle^{-\frac{11}{8}+\frac{\epsilon}{2}} +
\langle \tau \rangle^{-\frac{13}{8}+\frac{\epsilon}{2}}\right)\dtau\mathcal{F}^{2}(t)
\ls \langle t \rangle^{-\frac{5}{4}}\mathcal{F}^{2}(t).
\end{equation}
Note that
\[
\|\mathbf{u}\cdot\nabla\theta\|_{W^{8,1}}\lesssim \|u_{1}\partial_{1}\theta\|_{W^{8,1}} + \|u_{2}\partial_{2}\theta\|_{W^{8,1}}
\]
\[
\ls \|u_{1}\|_{L^{2}}\|\partial_{1}\theta\|_{H^{8}} + \|u_{1}\|_{H^{4}}\|\partial_{1}\theta\|_{L^{2}}
\]
\[
+ \|u_{2}\|_{L^{2}}\|\partial_{2}\theta\|_{H^{8}} + \|u_{2}\|_{H^{8}}\|\partial_{2}\theta\|_{L^{2}}
\]
\[
\ls \|u_{1}\|_{L^{2}}\|\partial_{1}\theta\|_{H^{\frac{6+2\delta}{\delta}}}^{\delta}\|\partial_{1}\theta\|_{H^{2}}^{1-\delta} + \|u_{1}\|_{H^{\frac{6+2\delta}{\delta}}}^{\delta}\|u_{1}\|_{H^{2}}^{1-\delta}\|\partial_{1}\theta\|_{L^{2}}
\]
\[
+ \|u_{2}\|_{L^{2}}\|\partial_{2}\theta\|_{H^{8}}^{\frac{1}{2}}\|\partial_{2}\theta\|_{L^{2}}^{\frac{1}{2}} + \|u_{2}\|_{H^{\frac{6+2\delta}{\delta}}}^{\delta}\|u_{2}\|_{H^{2}}^{1-\delta}\|\partial_{2}\theta\|_{L^{2}}.
\]
Taking $\delta$ such that $0<\delta  < \frac{1}{5}$ and $\frac{5\delta}{4}+2\delta\epsilon < \frac{1}{4}$, we obtain for $m>32$,
\[
O_{4}
\ls
\int_{0}^{t}\langle t-\tau\rangle ^{-\frac{5}{4}}\langle \tau\rangle^{-\frac{3}{2}+\frac{3\delta}{4}+\delta\epsilon}\dtau\mathcal{F}_{1}^{\delta}(t)\mathcal{F}_{3}^{2-\delta}(t)
\]
\[
+
\int_{0}^{t}\langle t-\tau\rangle ^{-\frac{5}{4}}\langle \tau\rangle^{-\frac{11}{8}+\frac{\epsilon}{2}}\dtau\mathcal{F}_{1}^{\frac{1}{2}}(t)\mathcal{F}_{3}^{\frac{3}{2}}(t)
+
\int_{0}^{t}\langle t-\tau\rangle ^{-\frac{5}{4}}\langle \tau\rangle^{-\frac{3}{2}+\frac{5\delta}{4}+\delta\epsilon}\dtau\mathcal{F}_{1}^{\delta}(t)\mathcal{F}_{3}^{2-\delta}(t)
\]
\begin{equation}\label{DJD4}
\ls \langle t\rangle^{-\frac{5}{4}}\mathcal{F}^{2}(t).
\end{equation}
Summing up (\ref{DJD1})-(\ref{DJD4}) yields
\begin{equation}\label{11theta}
\|\p_{11}\theta(t)\|_{L^{2}}=
\|\widehat{\p_{11}\theta}(t)\|_{\widehat{L}^{2}}\ls
\langle t\rangle^{-\frac{5}{4}}\left(\mathcal{F}_{0} +
\mathcal{F}_{0}^{2} + \mathcal{F}^{2}(t) +\mathcal{F}^{\frac{5}{2}}(t)+\mathcal{F}^{2+\delta}(t)\right).
\end{equation}

Next, we turn to $H^2$-estimate of $\omega$. According to (\ref{Eb6}), (\ref{JJJ2}) and (\ref{xb1}),
\[
\|\omega(t)\|_{H^{2}} = \|\widehat{\omega}(t)\|_{\widehat{H}^{2}}
\]
\[
\lesssim \left\|e^{-\nu\left(\xi^{2}+ \pi^2 k ^{2}\right) t}\widehat{\omega}_{0}\right\|_{\widehat{H}^{2}} +
\int_{0}^{t}\left\|e^{-\nu\left(\xi^{2}+ \pi^2k ^{2}\right)(t-\tau)}\widehat{\mathbf{u}\cdot\nabla\omega}(\tau)\right\|_{\widehat{H}^{2}} \dtau
\]
\[
+
\int_{0}^{t}\left\|e^{-\nu\left(\xi^{2}+ \pi^2k ^{2}\right)(t-\tau)}\widehat{\p_{1}\theta}(\tau)\right\|_{\widehat{H}^{2}}\dtau
\]
\[
\lesssim e^{-\nu t}\|\widehat{\omega}_{0}\|_{\widehat{H}^{2}} +
\int_{0}^{t}e^{-\nu(t-\tau)}\|\widehat{\mathbf{u}\cdot\nabla\omega}(\tau)\|_{\widehat{H}^{2}}\dtau +
\int_{0}^{t}e^{-\nu(t-\tau)}\langle\tau\rangle^{-\frac{3}{4}}\langle\tau\rangle^{\frac{3}{4}}\|\widehat{\p_{1}\theta}(\tau)\|_{\widehat{H}^{2}}\dtau
\]
\[
\lesssim e^{-\nu t}\|\omega_{0}\|_{H^{2}} +
\langle t\rangle^{-\frac{3}{4}}\left(\mathcal{F}_{0} +
\mathcal{F}_{0}^{2} +
\mathcal{F}^{2}(t)+\mathcal{F}^{\frac{5}{2}}(t)+\mathcal{F}^{2+\delta}(t)\right)
\]
\[
+\int_{0}^{t}e^{-\nu(t-\tau)}\|\left(\mathbf{u}\cdot\nabla\omega\right)(\tau)\|_{H^{2}}\dtau
\]
\begin{equation}\label{WWQ1}
\ls
\langle t\rangle^{-\frac{3}{4}}\left(\mathcal{F}_{0} +
\mathcal{F}_{0}^{2} +\mathcal{F}^{2}(t)+\mathcal{F}^{\frac{5}{2}}(t)+ \mathcal{F}^{2+\delta}(t)\right).
\end{equation}
Similarly, we use (\ref{11theta}) to obtain
\begin{equation}\label{WWQ3}
\|\p_{1}\omega(t)\|_{L^{2}}=\|\widehat{\p_{1}\omega}(t)\|_{\widehat{L}^{2}}
\ls \langle t\rangle^{-\frac{5}{4}}\left(\mathcal{F}_{0} +
\mathcal{F}_{0}^{2} +
\mathcal{F}^{2}(t)+
\mathcal{F}^{\frac{5}{2}}(t)+ \mathcal{F}^{2+\delta}(t)\right).
\end{equation}
Furthermore, applying (\ref{vel1})-(\ref{vel2}) and (\ref{WWQ1})-(\ref{WWQ3}) give
\[
\|u_{1}(t)\|_{H^{3}} = \left\|\frac{k\pi}{\xi^2+ \pi^2k^2}\widehat{\omega}(t)\right\|_{\widehat{H}^{3}}
\lesssim \|\widehat{\omega}(t)\|_{\widehat{H}^{2}}
\]
\[
\ls
\langle t\rangle^{-\frac{3}{4}}\left(\mathcal{F}_{0} +
\mathcal{F}_{0}^{2} +
\mathcal{F}^{2}(t)+\mathcal{F}^{\frac{5}{2}}(t)+ \mathcal{F}^{2+\delta}(t)\right),
\]
\[
\|u_{2}(t)\|_{H^{2}} = \left\|\frac{-i \xi}{\xi^2+ \pi^2k^2}\widehat{\omega}(t)\right\|_{\widehat{H}^{2}} \ls \|\widehat{\p_{1}\omega}(t)\|_{\widehat{L}^{2}}
\]
\[
\ls \langle t\rangle^{-\frac{5}{4}}\left(\mathcal{F}_{0} +
\mathcal{F}_{0}^{2} +\mathcal{F}^{2}(t)+\mathcal{F}^{\frac{5}{2}}(t)+ \mathcal{F}^{2+\delta}(t)\right),
\]
\[
\|\p_{1}u_{1}(t)\|_{H^{1}} \ls \|\widehat{\p_{1}\omega}(t)\|_{\widehat{L}^{2}}
\]
\[
\ls \langle t\rangle^{-\frac{5}{4}}\left(\mathcal{F}_{0} +
\mathcal{F}_{0}^{2} +
\mathcal{F}^{2}(t)+\mathcal{F}^{\frac{5}{2}}(t)+ \mathcal{F}^{2+\delta}(t)\right).
\]
This completes the proof of Lemma \ref{EL5}.
\hfill$\square$

\begin{Lemma}\label{EL6}
Let $m>32$. Then for $0<\delta <\frac{1}{5}$,
\[
\mathcal{F}_{4}(t) \lesssim \mathcal{F}_{0} + \mathcal{F}_{0}^{2} + \mathcal{F}^{2}(t) + \mathcal{F}^{\frac{5}{2}}(t)+ \mathcal{F}^{2+\delta}(t).
\]
\end{Lemma}
{\bf Proof.}
Using the fact that $\mathcal{L}_{2}(0) = 0$, we find from (\ref{Eb7}) that
\[
\partial_{t}\theta(x,y,t)= \partial_{t}\mathcal{L}_{1}(t)\theta_{0} + \partial_{t}\mathcal{L}_{2}(t)\left(\frac{1}{2}(-\Delta)\theta_{0} +\p_{t}\theta(x,y,0)\right)
\]
\begin{equation}\label{Eb41}
+ \int_{0}^{t}\partial_{\tau}\mathcal{L}_{2}(t-\tau)f_{2}(x,y,\tau)\dtau.
\end{equation}
From Lemma \ref{AL2}, one deduces  that
\[
\|\partial_{t}\theta(t)\|_{H^{3}}  =
\|\widehat{\partial_{t}\theta}(t)\|_{\widehat{H}^{3}}
\]
\[
 \lesssim
 \langle t\rangle^{-\frac{5}{4}}\|\omega_{0}\|_{W^{5,1}} +\langle t\rangle^{-\frac{5}{4}}\|\theta_{0}\|_{W^{8,1}} +
\langle t\rangle^{-\frac{5}{4}}\|\theta_{0}\|_{H^{m}}\|\omega_{0}\|_{H^{m}}
\]
\[
+ \int_{0}^{t}\langle t-\tau \rangle^{-\frac{5}{4}}\left(\|\left(\partial_{\tau}\mathbf{u}\cdot\nabla\theta\right)(\tau)\|_{W^{6,1}}+
\|\left(\mathbf{u}\cdot\nabla\partial_{\tau}\theta\right)(\tau)\|_{W^{6,1}}\right)\dtau
\]
\[
+\int_{0}^{t}\langle
t-\tau\rangle^{-\frac{5}{4}}\left(\|\partial_{1}(\mathbf{u}\cdot\nabla\omega)(\tau)\|_{W^{4,1}} +
\|\left(\mathbf{u}\cdot\nabla\theta\right)(\tau)\|_{W^{8,1}}\right)\dtau.
\]
Similar to $L^2$-estimate of $\partial_{11}\theta$ in Lemma \ref{EL5}, we find that for $0< \delta <\frac{1}{5}$,
\begin{equation}\label{xb7}
\|\partial_{t}\theta(t)\|_{H^{3}} =
\|\widehat{\partial_{t}\theta}(t)\|_{\widehat{H}^{3}} \lesssim
\langle t\rangle^{-\frac{5}{4}}\left(\mathcal{F}_{0} +
\mathcal{F}_{0}^{2} +
\mathcal{F}^{2}(t)
+\mathcal{F}^{\frac{5}{2}}(t)
+\mathcal{F}^{2+\delta}(t)\right).
\end{equation}

Taking time derivative on $(\ref{P4})_{1}$ and using Duhamel's principle,
\begin{equation}\label{vort1}
   \partial_{t}\omega(x,y,t)= e^{\nu t\Delta}\partial_{t}\omega(x,y,0) + \int_{0}^{t}e^{\nu(t-\tau)\Delta}\left(\partial_{1}\partial_{\tau}\theta - \partial_{\tau}(\mathbf{u}\cdot \nabla \omega)\right)(x,y,\tau)\dtau.
\end{equation}
Note that
\[
\partial_{t}\omega(x,y,0)= \nu\Delta\omega_{0} -\mathbf{u}_{0}\cdot \nabla \omega_{0} + \partial_{1}\theta_{0}.
\]
Moreover,
\[
\|\partial_{t}\omega(t)\|_{H^{2}} = \|\widehat{\partial_{t}\omega}(t)\|_{\widehat{H}^{2}}
\]
\[
\lesssim \left\|e^{-\nu\left(\xi^{2}+ \pi^2k^{2}\right) t}\widehat{\p_{t}\omega(0,x,y)}\right\|_{\widehat{H}^{2}} +
\int_{0}^{t}\left\|e^{-\nu\left(\xi^{2}+ \pi^2k^{2}\right)(t-\tau)}\widehat{\p_{\tau}(\mathbf{u}\cdot\nabla\omega)}(\tau)\right\|_{\widehat{H}^{2}} \dtau\]
\[
 +
\int_{0}^{t}\left\|e^{-\nu\left(\xi^{2}+ \pi^2
k^{2}\right)(t-\tau)}\widehat{\p_{1}\p_{\tau}\theta}(\tau)\right\|_{\widehat{H}^{2}}\dtau
\]
\[
\lesssim e^{-\nu t}\|\widehat{\p_{\tau}\omega(x,y,0)}\|_{\widehat{H}^{2}} +
\int_{0}^{t}e^{-\nu (t-\tau)}\|\widehat{\p_{\tau}(\mathbf{u}\cdot\nabla\omega)}(\tau)\|_{\widehat{H}^{2}}\dtau
\]
\[
+
\int_{0}^{t}e^{-\nu(t-\tau)}\langle\tau\rangle^{-\frac{5}{4}}\langle\tau\rangle^{\frac{5}{4}}\|\widehat{\p_{1}\p_{\tau}\theta}(\tau)\|_{\widehat{H}^{2}}\dtau
\]
\[
\lesssim e^{-\nu t}\left(\|\omega_{0}\|_{H^{4}} + \|\partial_{1}\theta_{0}\|_{H^{2}}+ \|\omega_{0}\|_{H^{3}}\|\mathbf{u}_{0}\|_{H^{3}}\right)
+\int_{0}^{t}e^{-\nu(t-\tau)}\|\p_{\tau}(\mathbf{u}\cdot\nabla\omega)(\tau)\|_{H^{2}}\dtau
\]
\begin{equation}\label{DDD1}
+
\langle t\rangle^{-\frac{5}{4}}\left(\mathcal{F}_{0} +\mathcal{F}_{0}^{2} +\mathcal{F}^{2}(t) + \mathcal{F}^{\frac{5}{2}}(t)
+\mathcal{F}^{2+\delta}(t)\right),
\end{equation}
where (\ref{xb7}) is used in the last inequality.
Now we need to estimate $\int_{0}^{t}e^{-\nu(t-\tau)}\|\p_{\tau}(\mathbf{u}\cdot\nabla\omega)(\tau)\|_{H^{2}}\dtau$. Obviously,
\[
\|\p_{t}(\mathbf{u}\cdot\nabla\omega)\|_{H^{2}}\ls \|\partial_{t}\mathbf{u}\cdot\nabla\omega\|_{H^{2}} + \|\mathbf{u}\cdot\nabla\partial_{t}\omega\|_{H^{2}}.
\]
By virtue of the interpolation inequality (\ref{BLT1}) in Lemma \ref{PL1},
\[
\|\partial_{t}\mathbf{u}\cdot\nabla\omega\|_{H^{2}}
\ls \|\partial_{t}\mathbf{u}\|_{H^{2}}\|\nabla\omega\|_{H^{2}}
\ls \|\partial_{t}\mathbf{u}\|_{H^{2}}\|\omega\|_{H^{6}}^{\frac{1}{2}}\|\omega\|_{L^{2}}^{\frac{1}{2}}.
\]
Then
\[
\int_{0}^{t}e^{-\nu(t-\tau)}\|\left(\p_{\tau}\mathbf{u}\cdot\nabla\omega\right)(\tau)\|_{H^{2}}\dtau
\]
\begin{equation}\label{DDD2}
\ls
\int_{0}^{t}
e^{-\nu(t-\tau)}\langle \tau\rangle^{-\frac{13}{8}+\frac{\epsilon}{2}}
\dtau\mathcal{F}_{4}(t)\mathcal{F}_{1}^{\frac{1}{2}}(t)\mathcal{F}_{3}^{\frac{1}{2}}(t)
\ls  \langle t\rangle^{-\frac{5}{4}} \left(\mathcal{F}^{2}(t)+\mathcal{F}^{\frac{5}{2}}(t)\right).
\end{equation}
Similarly,
\begin{equation}\label{DDD3}
\int_{0}^{t}e^{-\nu(t-\tau)}\|\left(\mathbf{u}\cdot\nabla\p_{\tau}\omega\right)(\tau)\|_{H^{2}}\dtau
\ls \langle t\rangle^{-\frac{5}{4}}\left(\mathcal{F}^{2}(t)+\mathcal{F}^{\frac{5}{2}}(t)\right).
\end{equation}
Hence, we deduce from (\ref{DDD1})-(\ref{DDD3}) that
\begin{equation}\label{DWd1}
\|\partial_{t}\omega(t)\|_{H^{2}} = \|\widehat{\partial_{t}\omega}(t)\|_{\widehat{H}^{2}} \lesssim
\langle t\rangle^{-\frac{5}{4}}\left(\mathcal{F}_{0} +
\mathcal{F}_{0}^{2} +
\mathcal{F}^{2}(t)
+ \mathcal{F}^{\frac{5}{2}}(t)
+ \mathcal{F}^{2+\delta}(t)\right).
\end{equation}
By using (\ref{vel1})-(\ref{vel2}) and (\ref{DWd1}),
\[
\|\partial_{t}\mathbf{u}(t)\|_{H^{3}} \lesssim \left\|\frac{1}{(\xi^{2}+ \pi^2k^{2})^{\frac{1}{2}}}\widehat{\partial_{t}\omega}(t)\right\|_{\widehat{H}^{3}}
\ls \|\widehat{\partial_{t}\omega}(t)\|_{\widehat{H}^{2}}
\]
\[
\lesssim
\langle t\rangle^{-\frac{5}{4}}\left(\mathcal{F}_{0} +
\mathcal{F}_{0}^{2} +
\mathcal{F}^{2}(t)
+ \mathcal{F}^{\frac{5}{2}}(t)
+ \mathcal{F}^{2+\delta}(t)\right).
\]
Thus, we finish the proof of Lemma \ref{EL6}.
\hfill$\square$

\subsection{Proof of Theorem \ref{Mrt1}}
According to Lemma \ref{EL1}-\ref{EL6}, there exists a constant $C_{0}\ge 1$ depending only on $\nu$ and $m$ such that
\begin{equation}\label{Ep36}
\mathcal{F}(t) \leq C_{0}\left(\mathcal{F}_{0}+\mathcal{F}_{0}^{2}\right)+ C_{0}\left(\mathcal{F}^{2}(t) + \left(\mathcal{F}^{3}(t) + \mathcal{F}^{4}(t)\right)^{\frac{1}{2}}+ \mathcal{F}^{\frac{5}{2}}(t)+ \mathcal{F}^{2+\delta}(t)\right)
\end{equation}
for $0<\delta<\frac{1}{5}$.
Assume that
\begin{equation}\label{assin1}
\mathcal{F}_{0} = \|\theta_{0}\|_{W^{8,1}} + \|\theta_{0}\|_{H^{m+1}} +\|\omega_{0}\|_{W^{5,1}} + \|\omega_{0}\|_{H^{m}}\leq \epsilon_{0}
\end{equation}
with $\epsilon_{0}\in(0,1)$  to be determined later.
Then by the definition of $F_i(t), i=1,2,3,4$, there exists $C_1\ge 2C_0$ such that
\begin{equation}\label{DDLL}
\mathcal{F}(0)= \mathcal{F}_{1}(0)+\mathcal{F}_{2}(0)+\mathcal{F}_{3}(0)+\mathcal{F}_{4}(0)
\leq C_1\mathcal{F}_{0}\leq C_{1}\epsilon_{0}.
\end{equation}
From Proposition \ref{tp1}, there exists $T^{*}\in (0,\infty]$ such that
\[
(\omega, \theta) \in C([0,T^*); H^{m})\times C([0,T^*);H^{m+1}).
\]
To complete the proof of Theorem \ref{Mrt1}, we only need to show $T^{*}=\infty$. To this end, we assume $T^*<\infty$ such that
\begin{equation}\label{T*}
\limsup\limits_{t\to T^*}\mathcal{F}_1(t)=\infty.
\end{equation}
Define
\begin{equation}\label{Ttmax1}
\tilde{T} \triangleq \max \{t\in [0,T^{*}):\mathcal{F}(t)\leq 4C_{1} \epsilon_{0}\} < T^*.
\end{equation}
By choosing $\epsilon_{0}$ such that $4C_{1} \epsilon_{0}\leq 1$ and using (\ref{Ep36}),
\begin{equation}\label{DDL2}
\mathcal{F}(t)\leq 2C_{0}\mathcal{F}_{0}+ 5C_{0}\mathcal{F}^{\frac{3}{2}}(t).
\end{equation}
By taking $\epsilon_{0}$ so small that $40 C_0 C_{1}^{\frac{1}{2}}\epsilon_{0}^{\frac{1}{2}}\leq 1$, we obtain from (\ref{DDLL}) that
\[
\mathcal{F}(t)\leq C_{1}\epsilon_{0}+  40 C_0 C_{1}^{\frac{3}{2}}\epsilon_{0}^{\frac{3}{2}}\leq 2C_{1} \epsilon_{0},\,\text{ for all}~t \in [0,\tilde{T}),
\]
which contradicts to (\ref{T*}) and the definition of $\tilde{T}$. This in turn implies that $T^{*} = \infty$.

\section{Conclusions}
It is shown in this paper the asymptotic behavior of solutions to Boussinesq equations (\ref{intp1}) near the specific stationary solution (\ref{steady0}). To this purpose, we formulate the momentum equation in terms of the vorticity, which enforces us to choose slip boundary condition for the velocity. Also in order to use spectral analysis, we choose functional spaces admitting high order compatibility conditions on the boundary. Since the anisotropic structure of the linearized equations is already reflected on the decay rates, we do not pursue this anisotropy on working functional spaces as well as the initial data.  Such an approach seems applicable to three dimensional setting and will be given in our subsequent paper \cite{LD2}. However, many problems are still left open. For example, we do not know how to handle the case of vanishing Dirichlet boundary conditions for the velocity. It is also interesting to show asymptotic stability in lower order spaces, even without any exact rate of convergence. Finally, it seems a challenging problem to show stability of the general steady solutions $\vartheta_s(y)$ such that $\vartheta_s'(y)>0$. We hope to solve these problems in the future work.

\centerline{\bf Acknowledgement}
The research is supported by NSFC under grants No. 12071211, 11771206.

\end{document}